\documentclass{article}
\usepackage{amssymb}
\usepackage{amsmath}
\usepackage{amscd}
\usepackage{amsthm}
\usepackage[all]{xy}

\def\text{\rm}

\newcommand\inte{\operatorname{int}}
\newcommand\End{\operatorname{End}}
\newcommand\Lie{\operatorname{Lie}}
\newcommand\alg{\operatorname{alg}}
\newcommand\Psm[2]{\psi_{\rm inv}^{#1}(#2)}

\newcommand\GR{\mathcal{G}}

\newcommand\ind{\operatorname{ind}}

\newcommand\CC{\mathbb C}

\newcommand\QQ{\mathbb Q}
\newcommand\RR{\mathbb R}
\newcommand\ZZ{\mathbb Z}

\newcommand\lgg{\mathfrak g}
\newcommand\pa{\partial}

\newcommand\Aut{\operatorname{Aut}}

\newcommand\KER{\operatorname{KER}}
\newcommand\COK{\operatorname{COK}}
\newcommand\Ker{\operatorname{Ker}}

\newcommand\wt[1]{\widetilde{#1}}
\newcommand\wh[1]{\widehat{#1}}
\def\g{\gamma}
\def\E{{\mathcal E}}

\def\g{\gamma}
\def\E{{\mathcal E}}

\def\Id{\operatorname{Id}}

\def\Hom{\operatorname{Hom}}

\def\a{\alpha}
\def\b{\beta}
\def\f{\varphi}
\def\s{\sigma}
\def\t{\tau}
\def\la{\lambda}
\def\F{\Phi}
\def\G{\Gamma}
\def\r{\rho}
\def\cC{{\mathcal C}}
\def\Ob{{\mathcal O}b}

\def\R{\mathbb R}
\def\C{\mathbb C}
\def\Lra{\Leftrightarrow}

\newtheorem{teo}{Theorem}[section]
\newtheorem{lem}[teo]{Lemma}
\newtheorem{prop}[teo]{Proposition}
\newtheorem{cor}[teo]{Corollary}

\theoremstyle{definition}
\newtheorem{dfn}[teo]{Definition}

\def\indg{\ind_\GR}
\def\LS{longitudinally smooth }

\newcommand\CI{{\mathcal C}^{\infty}}
\newcommand\CO{{\mathcal C}^{\infty}}
\newcommand\Cont{{\mathcal C}}
\newcommand\CIc{{\mathcal C}^{\infty}_{\text{c}}}

\newcommand\Diff{\operatorname{Diff}}
\newcommand\mP{\mathcal  P}

\begin{document}
\title{An index for gauge-invariant operators\\
and the Dixmier-Douady invariant}

\author{Victor Nistor\thanks{Partially supported by the NSF Young
Investigator Award DMS-9457859, NSF Grant DMS-9971951, and NSF Grant
DMS-9981251. }\ ~and Evgenij Troitsky\thanks{Partially supported by
RFFI Grant 99-01-01202 and Presidential Grant 00-15-99263.}}

\date{\today}


\maketitle

\begin{abstract}\ Let $\GR \to B$ be a bundle of compact Lie
groups acting on a fiber bundle $Y \to B$. In this paper we introduce
and study gauge-equivariant $K$-theory groups $K_\GR^i(Y)$.  These
groups satisfy the usual properties of the equivariant $K$-theory
groups, but also some new phenomena arise due to the topological
non-triviality of the bundle $\GR \to B$.  As an application, we
define a gauge-equivariant index for a family of elliptic operators
$(P_b)_{b \in B}$ invariant with respect to the action of $\GR \to B$,
which, in this approach, is an element of $K_\GR^0(B)$. We then give
another definition of the gauge-equivariant index as an element of
$K_0(C^*(\GR))$, the $K$-theory group of the Banach algebra
$C^*(\GR)$.  We prove that $K_0(C^*(\GR)) \simeq K^0_\GR(\GR)$ and
that the two definitions of the gauge-equivariant index are
equivalent.  The algebra $C^*(\GR)$ is the algebra of continuous
sections of a certain field of $C^*$-algebras with non-trivial
Dixmier-Douady invariant. The gauge-equivariant $K$-theory groups are
thus examples of twisted $K$-theory groups, which have recently proved
themselves useful in the study of Ramond-Ramond fields.
\end{abstract}

\markboth{\hfill \footnotesize \rm VICTOR NISTOR AND EVGUENI TROITSKY
\hfill}%
{\hfill \footnotesize \rm INDEX FOR GAUGE-INVARIANT OPERATORS \hfill}

\tableofcontents

\section*{Introduction}

Families of elliptic operators invariant with respect to a family of
Lie groups appear in gauge theory and in the analysis of geometric
operators on certain non-compact manifolds. They arise, for example,
in the analysis of the Dirac operator on a compact $S^1$-manifold $M$,
provided that we desingularize the action of $S^1$ by replacing the
original metric $g$ with $\phi^{-2}g$, where $\phi$ is the length of
the infinitesimal generator of the $S^1$-action. The main result of
\cite{NistorDOP} states that the kernel of the new Dirac operator on
the open manifold $M \smallsetminus M^{S^1}$ is naturally isomorphic
to the kernel of the original Dirac operator. A natural question to
ask then is when is this new Dirac operator (on $M \smallsetminus
M^{S^1}$) Fredholm.  The answer \cite{LMN,LN} is that, in general, the
Fredholm property of elliptic geometric operators on $M \smallsetminus
M^{S^1}$ is controlled by the invertibility of a certain family of
operators invariant with respect to a family of solvable Lie groups.
Operators invariant with respect to a family of Lie groups $\GR \to B$
will be called {\em gauge-invariant operators} or {\em $\GR$-invariant
operators} in what follows.

Gauge-invariant operators were studied in \cite{NistorFam}, where an
index theorem for these operators was also obtained in the case when
$\GR \to B$ is a family of solvable Lie groups. In this paper, we
study the case of $\GR$-invariant operators when $\GR$ is a bundle of
compact Lie groups.  In a certain way, the general case of operators
invariant with respect to a family of connected Lie groups can be
reduced to the cases when the fibers are solvable or compact groups,
by using the structure of connected Lie groups. There exist, however,
quite significant differences between families of solvable Lie groups
and families of compact Lie groups.

Here are two of the most important differences between families of
compact Lie groups and families of solvable Lie groups. First, a
family of compact Lie groups is always locally trivial.  In
particular, all fibers of a family of compact Lie groups are
isomorphic over each connected component of the base, see Section
\ref{S.Inv.Ops}. This need not be the case for families of nilpotent
groups, for example. Second, the theory of families of compact Lie
groups is more closely related to gauge theory and Physics than to
Analysis on non-compact manifolds. This is mainly because bundles of
compact Lie groups are one of the main objects of study in gauge
theory (see, for example, \cite{AtiyahBott, JaffeTaubes, WittenG,
DF}), but there are also other reasons.  For example, gauge-invariant
operators are related to anomalies in physics and to the Aronov-Bohm
effect \cite{Alvarez, ASZ, Belavin, BF, Brylinski, Moore,
Quillen}. Also, it is also possible that gauge-invariant operators are
related to the second quantization in field theory \cite{DF,
GlimmJaffe}.  An interesting possibility is to study whether
Ramond-Ramond charges \cite{FH, FW, MW, W3, W4} can be realized as
gauge-equivariant indices.

{\em From now on and throughout this paper, we shall assume that $\GR
\to B$ is a bundle of \underline{compact} Lie groups (except in
Section \ref{S.Inv.Ops} or when explicitly mentioned otherwise).}  The
relation between gauge-invariant operators and Physics provides
further motivation for the study of gauge-invariant operators.

A first problem that we solve for gauge-invariant operators is to
define their index. We actually give two definitions of this
index. The first definition is geometric and the second one is
algebraic. The advantage of the geometric definition is that it is
closer to the classical definition of the index of a family of
operators in the Atiyah and Singer paper on families \cite{AS4}. The
advantage of the algebraic definition of the index, however, is that
it works in general, whereas the geometric definition requires a
certain finite holonomy assumption on our bundle of Lie groups. This
finite holonomy assumption is automatically satisfied if the typical
fiber of our bundle of Lie groups has a center of dimension $\le 1$,
but not in general.

For the geometric definition of the gauge-equivariant index, we first
define groups $K^j_\GR(Y)$ for any fiber bundle $\pi: Y \to B$ on
which $\GR$ acts smoothly. (We assume that the action of $\GR$
preserves the fibers of $Y \to B$.)  These groups are defined
geometrically in terms of $\GR$-equivariant vector bundles on $Y$ and
give rise to a contravariant functor in $Y$ that has all the usual
properties of equivariant $K$-theory: homotopy invariance, continuity,
and Bott periodicity.  These groups behave well when $\GR$ has finite
holonomy (or, more generally, representation theoretic finite
holonomy, see Section \ref{Sec.Gauge.equiv} where these conditions
were introduced).  The geometric definition of the index associates an
element of the group $K^0_\GR(B)$ to any $\GR$-invariant continuous
family of pseudodifferential operators on $Y$. This is done as in the
classical case by perturbing our family to a family consisting of
operators whose range is closed (or, equivalently, such that the
family of their kernels defines a vector bundle over the base). The
difficulty is that we have to perform all these constructions
equivariantly.

For the algebraic construction of the analytic index of a
$\GR$-equivariant family of elliptic operators, we consider
$C^*(\GR)$, the $C^*$-algebra of $\GR$, which (in our case) can be
defined as the completion of the convolution algebra of $\GR$ with
respect to the action on each $L^2(\GR_b)$. Then
\begin{equation}\label{eq.isom.g.a}
        K^j_\GR(B) \cong K_j(C^*(\GR))
\end{equation}
whenever $\GR$ satisfies the finite holonomy assumption mentioned
above (Theorem \ref{thm.equiv.cat}).  Using some basic constructions
in $K$-theory, we can then give a direct definition of the
gauge-invariant index with values in $K_0(C^*(\GR))$, which turns out
to be equivalent to the geometric definition of the gauge-equivariant
index, if we take into account the isomorphism of Equation
\eqref{eq.isom.g.a} above. For this construction, we need no
assumption on $\GR$ (except that the fibers of $\GR \to B$ are compact
Lie groups), but it has the disadvantage that it is less elementary
and certainly less explicit. We leave the problem of determining the
gauge-equivariant index of a $\GR$-invariant family of elliptic
operators for another paper \cite{NT2}, since it requires several new
techniques.

The algebra $C^*(\GR)$ turns out to decompose naturally as a direct
sum of (algebras of continuous sections of) fields of finite
dimensional algebras on finite coverings of $B$. Assume that the
typical fiber $G$ of $\GR \to B$ is connected and denote by $G'$ the
derived group of $G$.  The Dixmier-Douady invariants of these fields
can be recovered {}from a unique class in $H^2(B, Z(G)\cap G')$ (see
also \cite{Dadarlat, Dixmier, DonovanKaroubi, Rosenberg3} for more on
these invariants). The $K$-theory groups of these algebras are twisted
$K$-theory groups, as the ones appearing in the study of Ramond-Ramond
fields \cite{Mathai, FH, FW, W3, W4}. This suggests some possible
connections between the structure of Ramond-Ramond fields and
gauge-equivariant index theory. An interesting possibility would be
that the Ramond-Ramond charges can be realized as gauge-equivariant
indices.

Let us now briefly describe the contents of each section. In Section
\ref{S.Inv.Ops}, we introduce bundles of Lie groups and define their
actions on spaces. We also introduce the algebras of (families of)
gauge-invariant pseudo-differential operators. This is the only
section were we do not assume that the fibers of $\GR \to B$ are
compact. In Section \ref{S.F.Hol}, we discuss the holonomy of the
representation spaces associated to bundles of Lie groups. For
instance, we introduce the condition that a bundle of Lie groups have
finite holonomy, which will play an important role in the study of
gauge-equivariant $K$-theory. The gauge-equivariant $K$-theory groups
$K^0_\GR(Y)$ of a bundle $Y \to B$ on which $\GR$ acts are introduced
in Section \ref{Sec.Gauge.equiv}. In that section, we establish
several properties of these groups. If the bundle $\GR \to B$ has
finite holonomy, then the groups $K^0_\GR$ behave like the usual
equivariant $K$-theory groups, but they can behave very differently if
this condition is not satisfied. For example, not every
$\GR$-equivariant vector bundle can be realized as a sub-bundle of a
trivial bundle in general; this is possible, however, if $\GR \to B$
has representation theoretic finite holonomy, which is actually the
reason why we introduced this condition, in the first place.  The
gauge-equivariant $K$-theory groups appear naturally as the receptacles
of the indices of gauge-equivariant families of elliptic
operators. These indices are defined geometrically in Section
\ref{Sec.geom} and analytically in Section \ref{S.An.Ind}. These two
definition are shown to coincide using the structure of the $K$-theory
of $C^*(G)$ and, especially, the isomorphism $K^0_\GR(B)=
K^0(C^*(\GR))$. In the Appendix we recall several basic constructions
in $K$-theory.

We would like to thank the Max Planck Institut for Mathematics in Bonn
for its kind support and hospitality while part of this work has been
completed. We would also like to thank M. Karoubi and A. Rosly for
useful discussions at Max Planck Institute.  Also, we would like to
thank M. Dadarlat for providing the example included in Section
\ref{Sec.structure.K}.

\section{Gauge-invariant pseudodifferential operators\label{S.Inv.Ops}}

We now describe the settings in which we shall work.  If $f_i : Y_i
\to B$, $i = 1, 2$, are two maps, we shall denote by
\begin{equation}
        Y_1 \times_B Y_2 := \{ (y_1,y_2) \in Y_1 \times Y_2,\,
        f_1(y_1) = f_2(y_2)\, \}
\end{equation}
their fibered product.

Fiber bundles will figure prominently in this paper. Even if we are
often interested primarily in smooth fiber bundles, we have found it
more convenient to work in the framework of ``longitudinally smooth''
fiber bundles, whose definition is recalled below. This simplifies and
makes more natural certain constructions in this paper, The reader can
safely assume for most results that the bundles are smooth. In fact,
we have included in ``$[\;\;]$'' the corresponding statements for
smooth fiber bundles.

\begin{dfn}\
Let $B$ be a locally compact topological space.  A locally trivial
fiber bundle $\pi: Y \to B$ with typical fiber $F$ is called {\em
longitudinally smooth\/} if, by definition, $F$ is a smooth manifold
and the structure group of this bundle reduces to $\Diff(F)$, the
group of $\CI$-diffeomorphisms of the fiber $F$.
\end{dfn}

In particular, every smooth bundle is longitudinally smooth. We now
introduce bundles of Lie groups.

\begin{dfn}
Let $B$ be a compact Hausdorff [respectively, smooth manifold], and
let $G$ be a Lie group. We shall denote by $\Aut(G)$ the group of
automorphisms of $G$. A {\em [smooth] bundle Lie groups $\GR$ with typical
fiber $G$ over $B$} is, by definition, a [smooth] fiber bundle $\GR
\to B$ with typical fiber $G$ and structure group $\Aut(G)$.
\end{dfn}

{\em We do not assume in this section that $\GR \to B$ has compact
fibers.}

We see thus that a [smooth] bundle of Lie groups with typical fiber
$G$ over $B$ is a longitudinally smooth [respectively, smooth] fiber
bundle over $B$ with typical fiber $G$ whose structure group reduces
from $\Diff(G)$ to $\Aut(G)$.

Let $\GR \to B$ and $\GR' \to B'$ be bundles of Lie groups. A {\em
morphism of bundles of Lie groups} (or simply, a {\em morphism}) $\g :
\GR \to \GR'$ is a continuous map covering a continuous map $\g_1 : B
\to B'$ such that the induced map $\GR_b \to \GR'_{\g_1(b)}$ is a
group morphism for any $b \in B$.

Let $d : \GR \to B$ be a bundle of Lie groups.  We shall often let
$\GR$ act on spaces $Y$. Let $\pi:Y \to B$ be a fiber bundle. We say
that $\GR$ {\em acts continuously on $Y$} if each group $\GR_b$ acts
continuously on $Y_b := \pi^{-1}(b)$ and the induced map $\mu$
$$
        \GR \times_B Y :=\{ (g,y) \in \GR \times Y, d(g) = \pi(y)\}
        \ni\, (g,y) \longrightarrow \mu(g,y) := gy \, \in Y,
$$
continuous. We shall also say that $Y$ is a {\em $\GR$-fiber bundle}.
If $\pi : Y \to B$ is longitudinally smooth bundle and $\GR_b \times
Y_b \to Y_b$ are smooth, then we shall say that $Y \to B$ is a {\em
longitudinally smooth $\GR$-fiber bundle}. If $\GR$ is a smooth bundle
of Lie groups, $Y \to B$ is a smooth fiber bundle and the map $\mu$ is
smooth, we shall say that $Y \to B$ is a {\em smooth $\GR$-fiber
bundle}.

Assume for the rest of this section that the quotient $Y /\GR := \cup
Y_b/\GR_b$ is compact.  Then we shall denote by $\Psm {m}Y$ the space
of continuous [respectively, smooth] families $D = (D_b)$, $b \in B$,
of order $m$, classical pseudodifferential operators acting on the
fibers of $Y \to B$ such that each $D_b$ is invariant with respect to
the action of the group $\GR_b$. Unless mentioned otherwise, we assume
that these operators act on half densities along each fiber. Then
$\Psm {\infty}Y := \cup_{m \in \ZZ} \Psm{m}Y$ is an algebra, by
classical results \cite{AS4, Hormander3}. Note also that $\Psm {m}Y$
also makes sense for $m$ not an integer.

We now discuss the principal symbols of operators in $\Psm
\infty{Y}$. Let
$$
        T_{vt}Y:=\ker (\pi_* : TY \to TB)
$$
be the bundle of {\em vertical} tangent vectors to $Y$, and let
$T_{vt}^*Y$ be its dual. We fix compatible metrics on $T_{vt}Y$ and
$T_{vt}^*Y$, and we define $S_{vt}^*Y$, the {\em cosphere bundle of
the vertical tangent bundle} to $Y$, to be the set of vectors of
length one of $T_{vt}^*Y$. Also, let
$$
        \sigma_m : \Psi^m(Y_b) \to \CI(S^* Y_b)
$$
be the usual principal symbol map, defined on the space of
pseudodifferential operators of order $m$ on $Y_b$ and with values
smooth functions on the unit cosphere bundle of $Y_b$. The definition
of $\sigma_m$ depends on the choice of a trivialization of the bundle
of homogeneous functions of order $m$ on $T^*_{vt}Y$, regarded as a
bundle over $S_{vt}^*Y$.  The principal symbols $\sigma_m(D_b)$ of an
element (or family) $D=(D_b) \in \Psm {m}Y$ then gives rise to a
smooth function on $\CO(S_{vt}^*Y)$, which is invariant with respect
to $\GR$, and hence descends to a smooth function on $S_{vt}^*Y$.  The
resulting function,
\begin{equation}
        \sigma_m(D) \in \Cont(S_{vt}^*Y))^\GR,
\end{equation}
will be referred to as the {\em principal symbol} of an element (or
operator) in $\Psm {m}Y$. Note that, $\sigma_m(D) \in
\CI(S_{vt}^*Y)^\GR$ when $\GR$ and $Y$ are smooth fiber bundles. As
usual, an operator $D \in \Psm {m}Y$ is called {\em elliptic} if, and
only if, its principal symbol is everywhere invertible.

In the particular case when $Y = \GR$ and $\GR$ is a smooth bundle,
$\Psm {\infty}\GR$ identifies with convolution operators on each fiber
$\GR_b$ that have compactly supported kernels, are smooth outside the
identity, and have only conormal singularities at the identity. In
particular, $\Psm {-\infty}{\GR} = \CIc(\GR)$, with the fiberwise
convolution product.

All constructions and definitions above extend to operators acting
between sections of a $\GR$-equivariant vector bundle $E \to Y$.
Recall that a bundle $E \to Y$ is $\GR$-equivariant if $\GR$ acts on
the total space of $E$, the projection $E \to Y$ is $\GR$ equivariant,
and the induced action $E_y \to E_{gy}$ between the fibers of $E \to
Y$ is linear. The space of order $m$, classical, $\GR$-invariant
pseudodifferential operators acting on sections of $E$ will be denoted
by $\Psm {m}{Y;E}$. The same construction then generalizes also to
$\GR$-invariant operators acting between sections of two
$\GR$-equivariant vector bundles $E_i \to Y$, $i = 0,1$. The
resulting space of order $m$, $\GR$-invariant pseudodifferential
operators acting between sections of $E_0$ and $E_1$ will be denoted
by $\Psm {m}{Y;E_0,E_1}$.

Let $T_{vt}Y$ denote the bundle of vertical tangent vectors to the
fibration $\pi: Y \to B$. Assume that there are given a
$\GR$-invariant metric on $T_{vt}Y$ and a $\GR$-equivariant bundle $W$
of modules over the Clifford algebras of $T_{vt}Y$. Then a typical
example of a family $D=(D_b) \in \Psm {\infty}Y$ is that of the family
of Dirac operators $D_b$ acting on the fibers $Y_b$ of $Y \to B$.
(Each $D_b$ acts on sections of $W\vert_{Y_b}$, the restriction of the
given Clifford module $W$ to that fiber.)

Later on, we shall need the following structure theorem for smooth
$\GR$-fiber bundles $Y$ (that is, when $\GR$, $Y$, and the action of $\GR$
on $Y$ are smooth). Recall that $\pi : Y \to B$ is the structural
projection map. Fix $b_0 \in B$ and let $G = \GR_{b_0} := d^{-1}(b_0)$
and $F = Y_{b_0} := \pi^{-1}(b_0)$. Let $\Aut(G,F)$ be the group of
automorphisms of the pair $(G,F)$, that is, $\Aut(G,F)$ consists of
pairs $(\alpha,\beta)$, where $\alpha$ is an automorphism of $G$ and
$\beta : F \to F$ is a diffeomorphism satisfying $\beta(gf) =
\alpha(g)\beta(f)$. Note that $\Aut(G,F)$ acts on both $G$ and
$F$. Moreover, the map (action) $G \times F \to F$ is
$\Aut(G,F)$-equivariant.

\begin{teo}\label{thm.fiber.Y}\
Let $\GR \to B$ be a smooth bundle of Lie groups acting smoothly on
the smooth fiber bundle $Y \to B$, whose typical fiber is denoted by
$F$. Then there exists a principal $\Aut(G,F)$-bundle $\mathcal Q \to
B$ such that
$$
        \GR \cong \mathcal Q
        \times_B G\,, \quad Y \cong \mathcal Q \times_B Y,
$$
and the induced map $\GR \times_B Y \to Y$ is obtained {}from the
$\Aut(G,F)$-equivariant map $G \times F \to F$.
\end{teo}

\begin{proof}\
Let $\mathcal Q$ be the set of triples $(b,\alpha,\beta)$ where $b \in
B$, $\alpha : G \to \GR_b$ is a group isomorphism, and $\beta : F \to
Y_b$ is an equivariant map, in the sense that $\beta(gf) = \alpha(g)
\beta(f)$ (as in the definition of the group $\Aut(G,F)$). The
projection map $p:\mathcal Q \to B$ is given by the projection onto
the first coordinate. By definition, $p^{-1}(b_0) = \Aut(G,F)$, so in
particular, this fiber is not empty.

It is clear that $\Aut(G,F)$ acts simply, transitively on each
non-empty fiber $p^{-1}(b)$. To complete the proof, all we need then
is to check that $\mathcal Q$ is locally trivial. To this end, we may
assume that $\GR = B \times G$.  The action of $\GR$ on $Y$ reduces
then to an action of $G$.  Since $G$ is compact, we can choose a
$G$-invariant metric on $Y$. The Levi-Civita connection will give then
rise to a $G$-equivariant diffeomorphism $Y_b \simeq F$ of the fibers
$Y_b$, for $b$ close to $b_0$. This proves the local triviality.
\end{proof}\smallskip

Similar ideas can be used to prove that a family of Lie groups with
connected compact fibers is locally trivial, and hence that it defines a
bundle of Lie groups. Let us first define the concept of a family of Lie
groups,
a concept that we consider only in the smooth category.

\begin{dfn}\label{def.fam.L.gr}\ Let $B$ be a smooth manifold. A {\em
family of Lie groups} is a submersion$d:\GR \to B$ such that each
$\GR_b := \pi^{-1}(b)$ is a Lie group and the induced map
\begin{equation}
        \GR \times_B \GR := \{ (g',g) \in \GR \times \GR, d(g') =
        d(g)\} \ni\, (g',g) \longrightarrow g'g^{-1}\, \in \GR
\end{equation}
is differentiable.
\end{dfn}

We remark that $B$ embeds naturally in $\GR$ as the space of units of
the groups $\GR_b$. If $\GR$ is a smooth family of Lie groups, then
$B$ is a smooth submanifold of $\GR$.  We also remark that a bundle of
Lie groups is a particular case of a continuous family groupoid (see
\cite{LMN,LN} for definitions) whose space of units identifies with
$B$. Let us denote by $d,r : \GR \to B$ the domain and, respectively,
the range maps of this groupoid. Then $d = r$, for our groupoid. If we
work in the differentiable category, this groupoid is a differentiable
groupoid.

For the rest of this section and then throughout the paper, we shall
assume that the fibers of the bundles (or families) of Lie groups that
we consider are compact.  The following result is probably already known.

\begin{teo}\label{thm.fam=bdle}\
Suppose that $\GR \to B$ is a family of Lie groups, that $B$ is
connected, and that all groups $\GR_b := d^{-1}(b)$ are compact and
connected. Then $\GR \to B$ is a locally trivial fiber bundle with
structure group $\Aut(G)$.
\end{teo}

\begin{proof}\
It is enough to prove that $\GR \to B$ is locally trivial. Fix $b_0
\in B$, and let $t_0 \in \GR_{b_0}$ be a topological generator of a
maximal torus of $\GR_{b_0}$ (recall that this means that the least
closed subgroup of $\GR_{b_0}$ containing $t_0$ is a maximal
torus). Choose a smooth local section $t$ of $\GR$, such that $t(b_0)
= t_0$. Let $H_b$ be the set of elements of $\GR_b$ commuting with
$t(b)$.

Let $\lgg = \cup \Lie \GR_b$ be the union of the Lie algebras of the
fibers of $\GR \to B$. Then $\lgg$ maps naturally to $B$ and can be
identified with the restriction to $B$ of the vertical tangent bundle
to $\GR \to B$.

Assume that in a small neighborhood of $b_0$ all the groups $\GR_b$
have maximal tori of the same dimension. Let $H_b^0$ be the connected
component of $H_b$. Then $H^0_b$ contains a maximal torus of $\GR_b$,
has the same dimension as a maximal torus of $\GR_{b_0}$, and hence
$H_b^0$ is a maximal torus in $\GR_b$, for any $b$ in a small
neighborhood of $b_0$. Then the root spaces associated to these
maximal tori vary continuously with $b$ and hence the root systems of
all the groups $\GR_b$ are isomorphic. This proves that the bundle
$\lgg$ of Lie algebras is locally constant, by Serre's theorem. Since
$\pi_1(\GR_b)$ is also constant, we also obtain that all the groups
$\GR_b$ are isomorphic.

Let us verify now the assumption that we made above, that is that in a
small neighborhood of $b_0$ all the groups $\GR_b$ have maximal tori
of the same dimension.

The Lie algebra $\Lie H_b$ is the kernel of $ad_{t_b}$.  Since the
operators $ad_{t_b}$ depend continuously on $b$ (on any trivialization
of $\lgg$ as a vector bundle), the dimension of the kernels of
$ad_{t_b}$ is $\le$ the dimension of $H_{b_0}$ in a small neighborhood
of $b_0$. Since $H_b$ contains a maximal torus of $\GR_b$, for any
$b$, we obtain that the set of points $b$ such that the dimension of a
maximal torus of $\GR_B$ is $\le n$ is an open subset of $B$. We can
approach then $b_0$ with a sequence of distinct points $b_n$ such that
$\GR_{b_n}$ have maximal tori of the same dimension $l$, for any $n$.

Then the maximal tori of $H_{b_n}^0$ also have dimension $l$.  Fix a
metric on $\lgg$ and let $X_1, \ldots, X_d \in \Lie \GR_{b_0}$ be a
basis of $\Lie \GR_{b_0}$, which we extend to a basis of $\Lie
H_{b_n}^0$, at least for $n$ large. Let $\epsilon >0$ be small, but
fixed. If we denote by $g_{ijn} := \exp(\epsilon X_j(b_n))\exp(\epsilon
X_i(b_n)) \exp(-\epsilon X_j(b_n)) \exp(-\epsilon X_i(b_n))$, then
\begin{equation}
	dist (g_{ijn}, 1 ) \to 0
\end{equation}
as $n \to \infty$, for any $i$ and $j$. Let $\delta > 0$ be arbitrary.
This proves that for $n$ large all the elements $\exp(\epsilon
X_j(b_n))$ are in a $\delta$ neighborhood of a torus of dimension
$l$. By letting $\delta \to 0$ and $n \to \infty$, we obtain that $l$
coincides with the dimension of $H_{b_0}$.
\end{proof}\smallskip

\section{Finite holonomy conditions\label{S.F.Hol}}

{\em From now on and throughout the paper all bundles of Lie groups
that we shall consider will have compact fibers. Also, unless
explicitly otherwise mentioned, $B$ will be a compact space and $\GR
\to B$ will be a bundle of compact Lie groups.}

We shall now take a closer look at the structure of bundles of Lie
groups whose typical fiber is a compact Lie group $G$.  Let $\Aut (G)$
be the group of automorphisms of $G$.  By definition, there exists
then a principal $\Aut(G)$-bundle $\mP \to B$ such that $\GR \cong \mP
\times_{\Aut(G)} G : = (\mP \times G)/\Aut(G)$.

Let $\widehat \GR$ be the (disjoint) union of the sets
$\widehat{\GR}_b$ of equivalence classes of irreducible
representations of the groups $\GR_b$. Using the natural action of
$\Aut(G)$ on $\wh G$, we can naturally identify $\widehat \GR$ with
$\mP \times_{\Aut(G)} \wh G$ as fiber bundles over $B$.

Let $\Aut_0(G)$ be the connected component of the identity in $\Aut
(G)$. The group $\Aut_0(G)$ will act trivially on the set $\wh G$,
because the later is discrete. Let $H_R := \Aut(G) / \Aut_0(G)$ and
$\mP_0 := \mP/ \Aut_0(G)$.  Then $\mP_0$ is an $H_R$-principal bundle
and $\widehat \GR \simeq \mP_0 \times_{H_R} \wh G$.

Assume now that $B$ is a path-connected, locally simply-connected
space and fix a point $b_0 \in B$.  Then the bundle $\mP_0$ is is
classified by a morphism $\pi_1(B,b_0) \to H_R$ because the structure
group $H_R$ of this principal bundle is discrete.

The space $\widehat \GR$ will be called {\em the representation space
of $\GR$}, the covering $\widehat \GR \to B$ will be called {\em the
representation covering associated to $\GR$}. Fix arbitrary a base
point of $B$. If $B$ is path-connected and locally simply-connected,
then the resulting morphism
\begin{equation}
        \rho : \pi_1(B,b_0) \to H_R := \Aut(G)/ \Aut_0(G),
\end{equation}
will be called {\em the holonomy of the representation covering of
$\GR$}.

For our further reasoning, we shall sometimes need the following
finite holonomy condition.

\begin{dfn}\label{def.c.f.h}\
We say that $\GR$ has {\em representation theoretic finite holonomy}
if every $\sigma \in \widehat{\GR}$ is contained in a compact-open
subset of $\widehat \GR$.
\end{dfn}

In the interesting cases, the above condition can be
reformulated as follows.

\begin{prop}\label{prop.r.f.h}\ Assume that $B$ is path-connected and
locally simply-connected. Then $\GR$ has {\em representation theoretic
finite holonomy} if, and only if, the following condition is
satisfied:\ for any irreducible representation $\sigma$ of $G$, the
set $\pi_1(B,b_0) \sigma \subset \widehat G$ is finite.
\end{prop}

\begin{proof}\ Since $\widehat{\GR}$ is a covering of
$B$, its compact-open subsets coincide with the finite unions of the
connected components of $\widehat{\GR}$ that are finite coverings of
$B$ (i.e. cover $B$ finitely-many times to one).  Let $B_\sigma$ be
the connected component of $\widehat{\GR}$ containig a given point
$\sigma \in \widehat{\GR}$. The typical fiber of $B_\sigma \to B$ is
$\pi_1(B,b_0)\sigma$. The result now follows.
\end{proof}\smallskip

The above condition ensuring representation theoretic finite holonomy
are difficult to check directly, so we shall also consider the
following closely related condition.

\begin{dfn}\label{def.f.h}\ Assume $B$ is smooth and connected.
We shall say that $\GR$ has {\em finite holonomy} if the image
$H_{\GR,b}$ of $\pi_1(B,b_0)$ in $H_R := \Aut(G)/\Aut_0(G)$ is finite.
\end{dfn}

These two ``finite holonomy'' conditions are related as follows.

\begin{teo}\label{thm.G=conn} \
Let $\GR \to B$ be a bundle of compact Lie groups over a
smooth connected manifold $B$. If $\GR$ has finite holonomy then it has
also representation theoretic finite holonomy.  If the fibers of $\GR
\to B$ are connected, then the converse is also true.
\end{teo}

\begin{proof}\ Let $\sigma$ be an irreducible representation of $G=\GR_{b_0}$.
Then $\pi_1(B,b_0)\sigma = H_{\GR}\sigma$ and hence
$\pi_1(B,b_0)\sigma$ is a finite set.

To prove the converse, let us write $G = G_0'Z_0$, where $G_0' \subset
G$ is the connected component of the subgroup $G'$ generated by
commutators and $Z_0 \subset G$ is the connected component of the
center of $G$. This is possible because $G$ is connected. Then $G_0'
\cap Z_0 = A$ is a finite subgroup and $G \simeq (G_0' \times Z_0)/A$,
where $A$ is embedded diagonally in the two subgroups.  Every
automorphism of $G$ maps each of $G_0'$ and $Z_0$ to itself. This
shows that the group of automorphisms of $G$ identifies with the
subgroup of automorphisms of $G_0' \times Z_0$ that map $A$ to
itself. Moreover, we have a canonical morphism $\Aut(G) \to
\Aut(Z_0)$. Since $\Aut(G_0')/\Aut_0(G_0')$ is a finite group, it is
enough then to show that the image of $\pi_1(B,b_0)$ in $\Aut(Z_0)$ is
finite whenever $\GR$ has representation theoretic finite holonomy.

The set $\widehat{S_0}$ of irreducible representations of $Z_0$
identifies with a lattice in the vector space $V= \Lie Z_0$. Choose a
basis $v_1, \ldots, v_n \in \widehat{S_0}$ of $V$. By assumption, the
sets $S_j = \pi_1(B,b_0)v_j$ are finite. The action of $\pi_1(B,b_0)$
then defines a group morphism from $\pi_1(B,b_0)$ to $\prod
\Aut(S_j)$, the product of the permutation groups of the sets
$S_j$. Let $K$ be the kernel of this morphism, then the morphism
$\pi_1(B,b_0) \to \Aut(Z_0)$ factors through $K$. This is enough to
complete the proof of the fact that representation theoretic finite
holonomy implies finite holonomy for a bundle of compact Lie groups
$\GR$ with connected fibers.
\end{proof}\smallskip

Let us note the following two consequence of the above proof that will
be useful in proving our results on the Dixmier-Douady invariants in
Section \ref{Sec.structure.K}. Denote by $G_{\inte} \subset \Aut(G)$
the subgroup of inner automorphisms of $G$.


\begin{prop}\label{prop.Z=1}\
Suppose $G$ is a compact, connected Lie group.  Then $\Aut_0(G)$
consists of inner automorphisms. If the center of $G$ has dimension
$\le 1$, then $\Aut(G)/\Aut_0(G)$ is a finite group.
\end{prop}

\begin{proof}\ We shall use the notation introduced in the proof of
Theorem \ref{thm.G=conn}. First, since $G_0'$ is semi-simple and connected,
the connected component $\Aut_0(G_0')$ of $\Aut(G_0')$ consists of inner
automorphisms. Since the image $G_{\inte} \subset \Aut(G)$ of $G$
acting by inner automorphisms is contained in $\Aut_0(G_0')$, we
obtain that $G_{\inte} = \Aut_0(G_0') = \Aut_0(G)$.

If $\dim Z_0 \le 1$, the group of automorphisms of $Z_0$ is finite
and this is enough to complete the proof.
\end{proof}\smallskip

Here is an example that the finite holonomy condition is not always satisfied
by a bundle of compact Lie groups $\GR \to B$. Let $A$ be the matrix
\begin{eqnarray*}
A = \left [
\begin{array}{cc}
3 & 2 \\ 4 & 3
\end{array}
\right ] .
\end{eqnarray*}
This matrix induces an automorphism $\alpha$ of the compact torus $T =
S^1 \times S^1$ by the formula
$$
        \alpha(z,w) = (z^3w^2, z^4w^3).
$$
Let us consider on the unit circle $S^1$ a bundle of tori $\GR_A$ with
fiber $T$ and holonomy $A$. This bundle can be realized as the
quotient of $\RR \times T$ by the equivalence relation $(t + n, z, w)
\equiv (t, \alpha^n(z,w))$, $n \in \ZZ$. The morphism $\ZZ \simeq
\pi_1(S^1) \to \Aut(T)$ sends then a generator of $\ZZ$ to the
automorphism $\alpha$. Clearly the range of this morphism is not
finite. The only irreducible representation $\sigma$ of $T$ with the
property that $\pi_1(S^1)\sigma$ is finite is the trivial
representation.

\section{Gauge equivariant $K$-theory\label{Sec.Gauge.equiv}}

In this section, we define the gauge equivariant $K$-theory groups for
spaces endowed with the action of a bundle of compact Lie groups
$\GR$. As we shall see, these are the right $K$-theory groups for our
index calculations. We formulate everything for the case of (locally)
compact Hausdorff topological spaces and continuous fiber bundles.
However, one can easily extend the following results to the smooth
setting (we usually include the necessary changes inside square
brackets).

Let $d: \GR \to B$ be a bundle of compact Lie groups over a compact
space $B$ and let $\pi: Y \to B$ be a fiber bundle. We assume, as in
the previous section, that $Y$ is a $\GR$-fiber bundle (that is, that
$\GR$ acts on $Y$). Let $E$ be a {\em finite-dimensional} vector
bundle $E \to Y$ equipped with an action of $\GR$. Such a vector
bundle will also be called a {\em $\GR$-equivariant vector
bundle}. This implies, in particular, that $E_b:=E\vert_{Y_b}$ is an
ordinary $\GR_b$-equivariant vector bundle over the $\GR_b$-space
$Y_b$, for any $b \in B$.

\begin{dfn}\label{def.morphism.b}\
Let $E \to Y$ be a $\GR$-equivariant vector bundle and let $E' \to Y'$
be a $\GR'$-equivariant vector bundle, for two bundles of Lie groups
$\GR \to B$ and $\GR' \to B'$. A {\em morphism} $(\g,\f):(\GR',E') \to
(\GR,E)$ is a pair of morphisms
$$
        (\g,\f),\quad \g:\GR' \to\GR, \quad  \f:E'      \to E,
$$
assumed to satisfy
$$
        \f(ge)=\g(g) \f(e),\qquad e\in E'_b,\quad g\in\GR'_b.
$$
(A map $\f$ with these properties will be called {\em $\g$-equivariant}.)

{\rm [}All the above maps are assumed to be smooth when working in the smooth
category.{\rm ]}
\end{dfn}

In particular, the maps $(\g,\f)$ in the above definition give rise
also to a map $B' \to B$ and to a $\g$-equivariant map $Y' \to Y$.

As usual, if $\psi : B' \to B$ is a continuous [respectively, smooth]
map, we define the {\em inverse image\/} $(\psi^*\GR,\psi^*E)$ of a
$\GR$-equivariant vector bundle $E \to Y$ by $\psi^*\GR = \GR \times_B
B'$ and $\psi^*E = E \times_B B'$. A particular case of this
definition is when $\psi$ is an embedding, when it gives the
definition of the restriction of a $\GR$-equivariant vector bundle $E$
to a closed subset $B'\subset B$ of the base of $\GR$, yielding a
$\GR_{B'}$-equivariant vector bundle. Usually $\GR$ will be fixed,
however.

The set of isomorphism classes of equivariant vector bundles
$E$ on $Y$, as above, will be denoted by $\E_\GR(Y)$. On this set we
introduce a monoid operation, denoted ``$+$,'' using the
direct sum of vector bundles. This defines a monoid structure on the
set $\E_\GR(Y)$.

\begin{dfn}\label{def.K}\
Let $\GR \to B$ be a bundle of compact Lie groups acting
on the fiber bundle $Y \to B$. Assume $Y$ to be compact.
The {\em $\GR$-equivariant $K$-theory group}
$K^0_{\GR}(Y)$ is defined as the group completion of the monoid
$\E_\GR(Y)$.

If $E \to Y$ is a $\GR$-equivariant vector bundle on $Y$, we
shall denote by $[E]$ its class in $K^0_\GR(Y)$. Thus $K^0_\GR(Y)$
consists of differences $[E]-[E_1]$.
\end{dfn}

The groups $K^0_\GR(Y)$ will  be called {\em gauge
equivariant} $K$-theory groups, when we do not need to specify $\GR$.
If $B$ is reduced to a point, then $\GR$ is a Lie group, and the
groups $K^0_\GR(Y)$ reduce to the usual equivariant $K$-groups.
More generally, this is true if $\GR \simeq B \times G$ is
a trivial bundle. The familiar functoriality properties of the usual
equivariant $K$-theory groups extend to the gauge equivariant
$K$-theory groups.

\begin{teo}\label{thm.funct}\
Assume that the bundle of Lie groups $\GR \to B$ acts on a fiber bundle $Y \to B$
and that, similarly, $\GR' \to B'$ acts on a fiber bundle $Y' \to B'$. Let $\g : \GR \to \GR'$
be a morphism of bundles of Lie groups and $f : Y \to Y'$ be a $\g$-equivariant
map. Then we obtain a natural morphism
\begin{equation}
        (\g,f)^* : K^0_{\GR'}(Y') \to K^0_{\GR}(Y)
\end{equation}
(denoted also $f^* : K^0_\GR(Y) \to K^0_\GR(Y')$ if $\gamma$ is the identity morphism).
\end{teo}

\begin{proof}\
The pull-back operation preserves the direct sum of equivariant
vector bundles, so it induces a morphism of monoids $\E_{\GR'}(Y')
\to \E_{\GR}(Y)$. The morphism $(\g,f)^* : K^0_{\GR'}(Y') \to K^0_{\GR}(Y)$
is the group completion of this morphism.
\end{proof}\smallskip

We now proceed to establish the main properties of the gauge equivariant
$K$-theory groups. Most of them are similar to those of the usual equivariant
$K$-theory groups, but there are also some striking differences.

First, let us observe that we can extend as usual the definition of
the gauge-equivariant groups to non-compact $\GR$-fiber bundles $Y$. Let $Y$
be a $\GR$-fiber bundle. We shall denote then by $Y^+ := Y \cup B$ the space
obtained from $Y$ by one-point compactifying each fiber.

\begin{dfn}
Assume that the typical fiber of the $\GR$-fiber bundle $Y
\to B$ is a locally compact space. We define then
$$
        K_{\GR}(Y) := \Ker \{K_\GR(Y^+) \to K_\GR(B)\}.
$$
\end{dfn}

Note that the above groups are ``compactly supported''
$K$-groups. These are the only $K$-theory groups that we shall
consider when working with non-compact manifolds.

We now discuss induction.  Let $\GR \subset \GR'$ be a sub-bundle of
the bundle of Lie groups $\GR' \to B$.  Also, let $Y$ be a
$\GR$-fiber bundle and $Y'$ be a $\GR'$-bundle. The fibered product over
$\GR$, namely $\times_{\GR}$, is defined as the quotient of
$\times_B$, the fibered product over $B$, by the action of $\GR$, its
action on $\GR'$ being by right translations, namely $\GR'
\times_{\GR} Y := (\GR' \times_B Y)/\GR$.  Then to any
$\GR$-equivariant vector bundle $E \to Y$ we can associate a
$\GR'$-equivariant vector bundle $E' = \iota(E)$ over $Y' := \GR'
\times_{\GR} Y$ by the formula $E' := \GR' \times_{\GR} E$.  This
operation gives rise to a morphism
\begin{equation}
        \iota : K^0_{\GR}(Y) \to K^0_{\GR'}(\GR' \times_{\GR} Y),
\end{equation}
called the {\em induction} morphism.

\begin{teo}\label{thm.ind}\
The induction morphism $\iota : K^0_{\GR}(Y) \to K^0_{\GR'}(\GR
\times_{\GR'} Y)$ is an isomorphism. The inverse is given by the
restriction to $Y \subset Y'$ and $\GR \subset \GR'$.
\end{teo}

\begin{proof}\
Note that $Y$ identifies with the image of $B \times_B Y$ (same $B$)
in $Y' := \GR' \times_{\GR} Y$. Moreover, the map $Y \to Y'$ is
$\GR$-equivariant. This shows that it makes sense to consider the
restriction map $r : K^0_{\GR'}(Y') \to K^0_{\GR}(Y)$. It follows from
the definition of $\iota(E)$, for a $\GR$-equivariant bundle $E \to
Y$, that $r \circ \iota$ is the identity.

The proof of the fact that $\iota \circ r$ is the identity follows
from the fact that any $\GR$-equivariant isomorphism $\beta : E
\vert_Y \to F \vert_Y$ of two $\GR'$-equivariant bundles extends
uniquely to a $\GR'$-equivariant isomorphism $E \to F$. We shall use
this as follows. Let $E$ be a $\GR'$-equivariant bundle. Then take $F
:= \iota \circ r(E)$. By what we have proved, we know that $r(E) = E
\vert_Y \simeq F \vert_Y = r(F)$ as $\GR$-equivariant bundles. Then $E
\simeq F$, as wanted.
\end{proof}\smallskip

The groups $K^0_{\GR}(Y)$ can be fairly small if the holonomy of $\GR$
is ``large.'' This is a new fenomenon, not encountered in the usual
equivariant $K$-theory. Here is an example.

\begin{prop}\label{prop.deg}\
Let $\GR_A$ be the bundle of two
dimensional tori over the circle $S^1$ introduced at the end of
Section \ref{S.Inv.Ops}, then $K^0_{\GR_A}(S^1) \simeq K^0(S^1)$.
\end{prop}

If $G$ were a compact group, the analogous statement would be that
$K^0_G(point) \simeq K^0(point)$, which is clearly true only if $G$ is
trivial. This kind of pathologies are ruled out by considering only
bundles of Lie groups with representation theoretic finite holonomy.

Let us denote by $C(\GR)$ the group of continuous sections of $\GR$,
that is, the group of continuous maps $\gamma : B \to \GR$. The group
$C(\GR)$ acts on $\G (E)$ according to the rule $(\gamma s)(y)
=\gamma(b)s((\gamma(b))^{-1} y),$ where $\gamma \in C(\GR)$, $s \in \G
(E) $, and $y \in Y_b$. We now define a continuous map $Av_{\GR} : \G
(E) \to \G^\GR (E)$, where $\G^\GR (E)$ is the space of
$C(\GR)$-invariant sections of $E$, by the formula
$$
        \big [ Av_{\GR} (s) \big ](x) := \int_{\GR_b}
        s(g^{-1}x) dg, \mbox{ if } x\in Y_b,
$$
the measure on $\GR_b$ being normalized to have total mass one. In the
following, by a $\GR$-invariant section of a $C(\GR)$-module we shall
understand a $C(\GR)$-invariant section.

\begin{lem}\label{prodoljsec}\
Let $X$ be a closed $\GR$-invariant sub-bundle of a compact
$\GR$-fiber bundle $Y$, and let $E \to X$ be a $\GR$-equivariant
bundle.  Let $s'$ be a $\GR$-invariant cross-section of the
restriction $E\vert_X \to X$.  Then $s'$ can be extended to a
$\GR$-invariant cross-section of $E$.
\end{lem}

\begin{proof}\
We proceed as in the classical case (cf. \cite{AK}). First extend the
section $s'$ to a section $s_1$ of $E$ over the whole of $Y$, not
necessarily equivariant.  The desired extension is obtained by setting
$s = Av_{\GR}(s_1)$.
\end{proof}\smallskip

We have the following analog of the corresponding classical result (see \cite{AK}
or \cite{KaKn}, for example).

\begin{teo}\label{epim}\
Let $E$ and $F$ be $\GR$-equivariant vector bundles over $X$, and $\a:
E \to F$ a $\GR$-equivariant morphism such that ${\a}_x: E_x \to F_x$
is an epimorphism for all points $x \in X$. Then there exists a
$\GR$-equivariant morphism $\b: F \to E$ such that $\a \b = \Id_F$.
\end{teo}

\begin{proof}\
First, we can define a possibly not equivariant morphism of bundles
$\wt\b:F\to E$ such that $\a\wt\b=\Id_F$ (see, e.g. \cite[Theorem
I.5.13]{KaKn}).  Then, let us take $\b:=Av_\GR(\wt\b)$, which we
define by regarding $\wt \b$ as an element of the $\GR$-equivariant
vector bundle $\Hom(F,E)$. Then
 $$
	\a\b = \a Av_\GR(\wt\b) = Av_\GR( \a \wt\b) = Av_\GR(\Id)=\Id.
 $$
\end{proof}\smallskip

A $\GR$-equivariant vector bundle $E \to Y$ on a $\GR$-fiber bundle $Y \to B$,
$Y$ compact, is called {\em trivial} if, by definition, there exists a $\GR$-equivariant
vector bundle $E' \to B$ such that $E$ is isomorphic to the pull-back of $E'$ to $Y$.
Thus $E \simeq Y \times_B E'$.

\begin{teo}\label{stabbund}\ Assume that $\GR \to B$ has representation
theoretic finite holonomy. Let $Y\to B$ be a compact $\GR$-fiber bundle and
$E \to Y$ be a $\GR$-equivariant vector bundle.
Then there exists a $\GR$-equivariant vector bundle $V \to B$
and a $\GR$-equivariant vector bundle $E'$ such  that
$Y \times_B V \cong E \bigoplus E'$.
\end{teo}

\begin{proof}\
Let $Y_b$ be the fiber of $Y \to B$ above some point $b \in B$.  Let
us recall how to embed $E_b := E\vert_{Y_b}$ into a trivial bundle,
for each $b \in B$. Choose sections $s_1,\ldots,s_n$ of $E\vert_{Y_b}$
such that they generate a finite dimensional $\GR_b$-invariant
subspace $V_b$ and they generate $\Gamma(E_b)$ as a
$C(Y_b)$-module. Then there exists a $\GR_b$-equivariant map $C(Y_b)
\times V_b \to \Gamma(E_b)$ that is surjective. The required embedding
into a trivial bundle is then obtained by an application of Theorem
\ref{epim} for $B$ reduced to $b$.

A set of sections $s_1,\ldots,s_n$ as above will be
called a {\em generating set of sections}.

Let $V_b$ be the representation
space of $\GR_b$ defined by these sections. Because $\GR$ has
representation theoretic finite holonomy, there exists a $\GR$-equivariant
vector bundle $W(b) \to B$ such that the fiber of this bundle at $b$ is
a representation of $\GR_b$ containing $V_b$. There will exist then a
$C(B)$-linear map $\Phi_b :  \Gamma(W(b)) \to \Gamma(E)$
and sections $\xi_j \in \Gamma(W(b))$ such that
$\Phi_b(\xi_j)\vert_{Y_b} = s_j$. By averaging
with respect to $\GR$ (using the map $Av_{\GR}$), we can assume that
$\Phi_b$ is $\GR$-equivariant.

Then $\Phi(\xi_j)$ will define by restriction a set of generating sections
of $E_{b'}$, for $b'$ in a neighborhood $U_b$ of $b$. Cover $B$ with finitely
many such neighborhoods $U_{b_j}$, and let $V \to B$ be the direct sum of all
the bundles $W(b_j)$ and $\Phi := \oplus \Phi_{b_j}$.

Our construction then gives a $\GR$-equivariant map
\begin{equation}
        1 \otimes \Phi : C(Y) \otimes_{C(B)} \Gamma(W) \to \Gamma(E)
\end{equation}
that is surjective by construction. We have thus constructed a surjective
map $Y \times_B V \to E$. An application of Theorem \ref{epim}
concludes the proof.
\end{proof}\smallskip

Let us observe that the above result is not true for $\GR_A$, the bundle of two
dimensional tori over $S^1$ considered at the end of Section \ref{S.Inv.Ops}.
Indeed, let $\GR \subset \GR_A$ be subset of all elements of order two of the
fibers of $\GR_A$. Then $\GR \to S^1$ is a trivial bundle of finite groups:\
$\GR = S^1 \times A$, with $A \simeq (\ZZ/2\ZZ)^2$. Let
$$
        Y' := \GR_A \times_{\GR} S^1 = \GR_A/\GR.
$$
We know by Theorem \ref{thm.ind} that
$$
        K^0_{\GR_A}(Y') \simeq K^0_{A}(S^1) \simeq R(A) \otimes K^0(S^1) = R(A),
$$
the isomorphism being given by restriction to $S^1 \subset Y'$. On the other hand,
if $E \to Y'$ were a sub-bundle of a trivial $E'$ bundle over $Y'$,
then $E \vert_{S^1}$ would also be a $\GR$-equivariant
sub-bundle of the trivial $\GR$-equivariant bundle $E'\vert_{S^1}$. If $E''$ is a
$\GR_A$-equivariant bundle over $S^1$, then the pull-back to $Y'$ followed by the
restriction to $S^1$ corresponds to restricting the action of $\GR_A$ to an
action of $\GR$. Thus any bundle of the form $E'\vert_{S^1}$, with $E'$ a trivial
$\GR_A$-bundle, will be trivial over $S^1$ and will have the trivial action of $A$, by
Proposition \ref{prop.deg}. Any sub-bundle
of $E'$ will again have the trivial action of $A$. This shows that the $\GR_A$-equivariant
bundles over $Y'$  that can be realized as sub-bundles of trivial bundles have a class
in $K^0_{\GR_A}(Y') \simeq R(A)$ corresponding to multiples of the trivial representation.
This gives the following result.

\begin{prop}\label{prop.note}\
Thus not every $\GR_A$-equivariant bundle over $Y'$ can be realized as a sub-bundle of
a trivial bundle.
\end{prop}

We now check that the category of $\GR$-equivariant vector bundles is
a Banach category (see Definition~\ref{kar3}). All the other properties of the gauge
equivariant $K$-theory groups that we shall prove will turn out to be consequences
of a some general theorems on Banach categories from \cite{KarCli,KaKn}. We shall
obtain, in particular, that gauge equivariant $K$-theory
has long exact sequences and satisfies Bott periodicity.

\begin{prop}\label{prop.Banach}\ The category of $\GR$-equivariant vector
bundles over a $\GR$-fiber bundle $Y \to B$, $Y$ compact, is a Banach category.
\end{prop}

\begin{proof}\
First, the set $\G(E)$ of all continuous sections $s : Y \to E$ of a
$\GR$-equivariant vector bundle $E \to Y$ over a compact $\GR$-fiber bundle $Y \to B$
becomes a Banach space when endowed with the ``sup''-norm.
Consider now two $\GR$-equivariant vector bundles $E$ and $F$ over a compact
$\GR$-fiber bundle $Y$. The vector bundle $\Hom(E,F)$ will have a natural $\GR$-action.
As usual, we can identify  $\G^\GR (\Hom (E, F))$ with the set
of $\GR$-equivariant morphisms $\f: E \to F$, which is hence a Banach space.
The composition of morphisms $\Hom(E_1,E_2) \times \Hom(E_2,E_3) \to \Hom(E_1,E_3)$
is continuous because the category of vector bundles is a Banach category. The
restriction to $\GR$-equivariant morphisms will also be continuous. This
checks all conditions of Definition~\ref{kar3} (see the Appendix, where
Banach categories are discussed following \cite{KarCli,KaKn})), and hence the proof is now
complete.
\end{proof}\smallskip

This proposition allows now to establish several useful Lemmata.

\begin{lem}\label{prodoljmor}\
Let $E,F \to X$ be $\GR$-equivariant vector bundles and let $Y \subset X$
be a $\GR$-equivariant sub-bundle over $B$.
Let $\f': E\vert_Y \to F\vert_Y$ be a morphism of the
restrictions of the $\GR$-fiber bundles $E$ and $F$ to $Y$. Then $\f'$ can be extended
to a $\GR$-equivariant morphism $\f: E \to F$. If $\f$ is an isomorphism, then
there exists a $\GR$-invariant open neighborhood $U$ of $Y$
such that
$$
        \f\vert_U: E\vert_U \to F\vert_U
$$
is an isomorphism. Any two such extensions are homotopic to each other
in the class of $\GR$-equivariant isomorphisms over some
$\GR$-invariant neighborhood of $X$.
\end{lem}

\begin{proof}\
We apply Lemma~\ref{prodoljsec} to the bundle $\Hom(E,F)$ and use the
fact that the set of isomorphisms forms an open subset of the set of
all homomorphisms.
\end{proof}\smallskip

If $Y$ is a $\GR$-fiber bundle over $B$ and $I$ is the unit interval, then
we define the action of $\GR_b$ on $Y_b \times I$ by the formula
$g (y, t) = (gy, t)$. So we obtain a $\GR$-fiber bundle $Y\times I$ over
$B$. This is a particular case (for $Z=B\times I$ and the trivial
action of $\GR$ on $Z$) of the following construction.
Suppose that  $Y$ and $Z$ are $\GR$-fiber bundles over $B$. Then the bundle
$Y\times_B Z\to B$ can be equipped with the diagonal action of $\GR$.

\begin{lem}\label{homot}
Let $\pi_Y : Y \to B$ and $\pi_X : X \to B$ be compact $\GR$-fiber
bundles and $f_t: Y \to X$ be a continuous homotopy of
$\GR$-equivariant mappings $(0\le t\le 1)$ satisfying $\pi_X \circ f_t
= f_t \circ \pi_Y$. Suppose that $E$ is a $\GR$-equivariant vector
bundle over $X$. Then $f^*_0 E \cong f^*_1 E$.
\end{lem}

\begin{proof}\
We proceed as in \cite{AK}, Lemma 1.4.3. Denote by $I$ the unit
interval $[0,1]$.  Let $\pi : Y \times I \to Y$ be the first
projection, and let $f : Y \times I \to X$, $f(y,t) := f_t(y)$ be the
given homotopy.  Let us consider the bundles $f^*E$ and $\pi^* f^*_t
E$ for some fixed $t$.  Over $Y\times \{t\}$ we have a natural
isomorphism (identification) of these bundles. By
Lemma~\ref{prodoljmor}, since $Y$ is compact, there exists a
neighborhood $U_t$ of $t$ in $I$ such that $f^*E$ and $\pi^* f^*_t E$
are isomorphic over $U_t$. Hence the isomorphism class of
$(f^*E)\vert_t=f^*_t E$ is a locally constant function of $t$. The
connectedness of $I$ completes then the proof.
\end{proof}\smallskip

\begin{cor}\label{cor.herm}
Let $Y \to B$ be a $\GR$-fiber bundle and $E \to Y$ be a
$\GR$-equivariant vector bundle. Then there exists a $\GR$-invariant
Hermitian metric on $E$.
\end{cor}

\begin{proof}\
Choose an arbitrary Hermitian metric on $E$, regarded as an element of
the $\GR$-equivariant vector bundle $E^* \otimes E^*$. Its average
will then be a $\GR$-invariant Hermitian metric.
\end{proof}\smallskip

Another immediate consequence of Theorem~\ref{stabbund} is the
following important property.

\begin{teo}\label{quasisurj}\
Let $\r^{Y,Y'}: \E _\GR(Y) \to \E _{\GR}(Y')$ be the restriction
functor, where $Y'$ is a closed $\GR$-invariant sub-bundle of
$Y$. Then $\r^{Y,Y'}$ is a full quasi-surjective Banach functor in the
sense of {\rm~\ref {kar4}} .
\end{teo}

\begin{proof}\
We shall use the concepts recalled in the Appendix. Since the
restriction map
$$
	\r_*:\G^\GR(\Hom(E,F)) \to \G^\GR(\Hom(E|_{Y'},F|_{Y'}))
$$
is a continuous linear map of Banach spaces (cf. the proof of
Proposition~\ref{prop.Banach}), $\r$ is a Banach functor.  By Lemma
\ref{prodoljmor}, $\r$ is full. Finally, it is quasi-surjective by
Theorem~\ref{stabbund}, because $Y\times_B V$ extends $Y'\times_B V$.
\end{proof}\smallskip

\begin{dfn}\
\rm Following the general scheme presented in the Appendix (Section
\ref{banachcath}), we define the $K$-groups for a compact space $B$ by
setting
$$
        \begin{array}{rl} K^{p,q}_\GR(Y)&=K^{p,q}( \E_\GR(Y),\\[5pt]
        K^{n}_\GR(Y)&=K^{n}( \E _\GR(Y)),\\[5pt]
        K_\GR(Y)&=K^{0,0}(\E_\GR(Y)).  \end{array}
$$
\end{dfn}

By Theorem \ref{quasisurj}, setting $K^n_\GR(Y,Y') := K^n(\r^{Y,Y'})$
in the sense of Definition \ref{kar9a}, we get (by
\cite[II.3.22]{KaKn} and \cite[2.3.1]{KarCli}) a long exact sequence
\begin{multline}\label{exact}
        \dots \to K^{n-1}_\GR(Y,Y')\to K^{n-1}_\GR(Y)\to
        K^{n-1}_{\GR_1}(Y') \\ \to K^n_\GR(Y,Y') \to K^n_\GR(Y) \to
        K^n_{\GR}(Y') \to \dots
\end{multline}

\begin{teo}\label{bott} {\rm (Bott-Clifford periodicity) \/}
We have a natural isomorphism
$$
     K^n_\GR(Y,Y') \cong K^{n-2}_\GR(Y,Y').
$$
 \end{teo}

\begin{proof}\
The category of $\GR$-equivariant vector bundles is a Banach category
(see Definitions \ref{kar3} and \ref{karpsab} of the Appendix) This
and the periodicity of the Clifford algebras directly implies our
result (see, for example, \cite[\S III.4]{KaKn}).
\end{proof}\smallskip

\begin{teo} {\rm (Periodicity.) \/}
Let $Y \to B$ be a compact $\GR$-fiber bundle. We have a natural
isomorphism
$$
     K^1_\GR(Y) \cong K_\GR(Y\times D^1,Y\times S^0),
$$
where $(D^n, S^{n-1})$ is a ball and its boundary.
\end{teo}

\begin{proof}\
This is proved for general (complex) Banach categories in
\cite[Theorem 2.3.3]{KarCli}.
\end{proof}\smallskip

\begin{teo}\label{teo:excision}
Suppose that $Y\to B$ is a closed invariant sub-bundle of a compact
$\GR$-fiber bundle $X\to B$. Then the projection $\kappa:X\to X/Y$
induces an isomorphism $\kappa^*: K_\GR(X/Y,\{y\})\to K_\GR(X,Y)$.
\end{teo}

\begin{proof}\
We shall follow the proof in \cite[II.2.35]{KaKn}.

Let us prove first that $\kappa^*$ is surjective. Let $d(E,F,\a)$ be
an element of $K_\GR(X,Y)$.  By adding the same bundle to $E$ and $F$,
we may assume, without loss of generality, that $F\cong X\times_B V$
(that is, that $F$ is the pull-back from $B$ of a vector bundle $V \to
B$). We want to find a triple $(E',F',\a')$ defining an element of
$K_\GR(X/Y,\{y\})$, such that the triples
$(\kappa^*(E'),\kappa^*(F'),\kappa^*_1(\a'))$ and $(E,F,\a)$ are
isomorphic, where $\kappa_1:Y\to \{y\}$ maps $Y$ to a point.
According to Lemma~\ref{prodoljmor}, there is a closed $\GR$-invariant
neighborhood $N$ of $Y$ and an isomorphism $\b:E\vert_N\to F\vert_N$
such that $\b\vert_Y=\a$.  Let $E'$ be the vector bundle over $X/Y$
obtained by clutching the bundle $E\vert_{X\setminus Y}$ and the
bundle $n/Y\times_B V$, using $\b\vert_{X\setminus Y}$.  One has
$X\setminus Y = (X/Y) \setminus \{y\}$ and $N \setminus Y = (N/Y)
\setminus \{y\}$. Let $F' = X/Y \times_B V$, and let
$\a':E'\vert_{\{y\}} \to F'\vert_{\{y\}}$ be the isomorphism induced
by the above clutching.

Then we can define an isomorphism $f: E\to \kappa^*(E')$ by
$f\vert_{X\setminus Y}=\Id$, with the identification $\kappa^*
E'\vert_{X\setminus Y}= E'\vert_{X\setminus Y}=E\vert_{X\setminus Y}$,
and $f\vert_N=\Id$ with the identification $\kappa^*(E')\vert_N=
(X\times_B V)\vert_N$. It is now obvious that the diagram
$$
        \xymatrix{E\vert_Y \ar[r]^\a \ar[d]_{f\vert_Y} &
        F\vert_Y\ar@{=}[d]\\ \kappa^*(E')\vert_Y
        \ar[r]^{\kappa^*_1(\a')} &\kappa^*(F')\vert_Y}
$$
is commutative.

We now prove that $\kappa^*$ is injective.  Let $d(E',F',\a')$ be an
element of $K_\GR(X/Y,\{y\})$ such that
$$
        \kappa^*(d(E',F',\a')) =
        d(\kappa^*(E'),\kappa^*(F'),\kappa^*_1(\a'))=0.
$$
According to Proposition~\ref{karopisnulfun}, there is a bundle $T$
over $X$ such that $\kappa^*_1(\a')\oplus\Id\vert_{T\vert_Y}$ can be
extended by an isomorphism $\b:\kappa^*(E')\oplus T\to
\kappa^*(F')\oplus T$.  As before we may assume that $T=X\times_B
V$. Let $T'= X/Y\times_B V$. Let $\b':E'\oplus T'\to F'\oplus T'$ be
the isomorphism which is equal to $\b$ over $X\setminus Y$, and to
$\a'$ over $\{y\}$. Then $\b'$ is continuous and is an extension of
$\a'\oplus \Id_{T\vert_{\{y\}}}$ over $X/Y$. By
Proposition~\ref{karopisnulfun}, $d(E',F',\a')=d(E'\oplus T', F'\oplus
T', \a'\oplus \Id_{T'})=0$.
\end{proof}\smallskip

\begin{lem}\ Let $Y \to B$ be a compact $\GR$-fiber bundle.  Then
      $
     K_{\GR}(Y) \cong K_\GR(Y^+,B).
      $
\end{lem}

\begin{proof}\
It follows by writing the long exact sequence (Equation \eqref{exact})
for the pair $(Y^+,B)$.
\end{proof}\smallskip

\smallskip
Using the above Lemma, we see that the long exact sequence extends of
Equation \eqref{exact} extends to non-compact $\GR$-fiber bundles $X$
and $Y$.

\begin{dfn}\ \rm Assume $Y \to B$ is a possibly non-compact
$\GR$-fiber bundle. Let $Y^+ := Y \cup B$ be the fiberwise one-point
compactification. Let
      $$
K^{-1}_{\GR}(Y)=\Ker\{K^{-1}_\GR(Y^+) \to K^{-1}_\GR(B)\},
      $$
      $$
K^n_{\GR}(Y,Y')=K_{\GR}((Y \setminus Y') \times \R^n).
      $$
\end{dfn}

It is necessary to verify the compatibility condition
$K^{-1}_{\GR}(Y) \cong K_{\GR}(Y \times \R).$ For this purpose,
let us consider the bundle $Z:=Y\times \R_+$, where $\R_+=[0,+\infty)$
(c.f., \cite[Theorem II.4.8]{KaKn}).  Then $Z$ is fiberwise
homeomorphic to $Y^+\times [0,1]\setminus Y^+ \vee [0,1]$, where $1$
is the base point of $[0,1]$.  Hence $Z^+$ is fiberwise homeomorphic
to the fiberwise quotient $Y^+_b \times [0,1] / Y^+_b \vee [0,1]$.
Since $ Y^+_b \vee [0,1]$ is invariant, the identification is
equivariant.  Under this identification, let us define an equivariant
fiberwise homotopy $r:Z^+\times [0,1]\to Z^+$ by the formula
$r([y,t],u)=[y,1+(1-t)u]$, where $y\in Y^+$, $t,u \in [0,1]$, and
$[\:,\:]$ means the class in the quotient space.  When $u=0$, the
image is $B$. Hence $K_\GR(Y\times \R_+)=K^{-1}_\GR(Y\times \R_+)=0
$. The long exact sequence of the pair $(Y\times \R_+, Y)$
$$
        K^{-1}_\GR(Y\times \R_+)\to K^{-1}_\GR(Y)\to K_\GR(Y\times \R)
        \to K_\GR(Y\times \R_+),
$$
proves immediately the required compatibility:\ $K^{-1}_{\GR}(Y)
\cong K_{\GR}(Y \times \R).$ In this proof we have also used
the identifications
$$
        K_\GR(Y\times \R_+, Y)\cong K_\GR(Y\times \R_+\setminus
        Y)\cong K_\GR(Y\times \R),
$$
in view of Theorem \ref{teo:excision}.

 \section{The analytic index: a geometric approach\label{Sec.geom}}

For a family of elliptic operators invariant with respect to a bundle
of compact Lie groups $\GR \to B$, it is possible to extend the
definition of the family index to obtain an index with values in the
gauge-equivariant $K$-theory groups $K^0_{\GR}(B)$ introduced in the
previous section. In this section we provide an explicit geometric
construction of this index when $\GR$ has representation theoretic
finite holonomy. The general case requires different methods and will
be treated in Section \ref{S.An.Ind}. We continue to assume that $B$
is compact.

\begin{dfn}\
A {\em locally trivial bundle of Hilbert spaces} over $B$ is a fiber
bundle $H \to B$ such that $B$ can be covered with open sets
$U_\alpha$ with the property that $H\vert_{U_\alpha} \simeq U_\alpha
\times H_0$, for some fixed Hilbert space $H_0$ and the transition
functions are continuous in norm.

Let $\GR \to B$ be a bundle of Lie groups.
A {\em locally trivial bundle of $\GR$-Hilbert spaces} over $B$ is a
{\em locally trivial bundle of Hilbert spaces} $H \to B$ together with
a fiber-preserving action of $\GR$ on $H$ that consists of continuous
families in any trivialization of $H$.
\end{dfn}

It is known that every locally trivial bundle of {\em infinite
dimensional} Hilbert spaces on a finite dimensional base is actually
trivial, because the space of unitary operators of an infinite
dimensional Hilbert space is contractable (Kuiper's theorem
\cite{Kuiper}). See also Dixmier's book \cite{Dixmier}.

{\em We fix in this section a bundle of compact Lie groups $\GR \to
B$, with representation theoretic finite holonomy and with a $B$
compact, path connected, and locally simply-connected topological
space.}

\begin{lem}\label{lemma.KER}\ Assume that $\GR \to B$ has representation theoretic
finite holonomy (as above).
Suppose that $H^0\to B$ and $H^1\to B$ are two locally trivial bundles
of $\GR$-Hilbert spaces. Suppose also that $F=(F_b :H^0_b \to H^1_b)_{b \in
B}$ is a family of $\GR$-invariant Fredholm operators that is norm-continuous in
any trivialization of $H^i \to B$. Then there exists a finite-dimensional
$\GR$-invariant vector sub-bundle $\KER \subset H^0$ such that:
\begin{enumerate}
\item $F_b:(\KER_b)^\bot \to F_b((\KER_b)^\bot)$ is a $\GR_b$-isomorphism
for every $b\in B$;

\item $\COK:=\bigcup\limits_{b\in B}(F_b(\KER_b))^\bot\subset H^1$ is
a finite-dimensional $\GR$-invariant sub-bundle.
\end{enumerate}
\end{lem}

\begin{proof}\
For any irreducible representation $\sigma$ of $G$, let us denote
$S_\sigma = \pi_1(B,b_0)\sigma$. Our assumption that $\GR$ has
representation theoretic finite holonomy is equivalent to saying that
all sets $S_\sigma$ are finite (Proposition \ref{prop.r.f.h}).
These sets then partition $\widehat G$,
the set of irreducible representations of $G$.

For each $b \in \GR_b$ we obtain a subset $S_{\sigma,b} \subset
\widehat {\GR_b}$, defined as the fiber above $b$ of the space $\mP
\times_{\Aut(G)} S_\sigma$. Let us identify $G$ with one of the fibers
$\GR_{b_0}$ of $\GR \to B$. In particular, this identifies $S_\sigma$
with $S_{\sigma,b_0}$.  Since $B$ is path connected, the set
$S_{\sigma,b}$ of irreducible representations of $\GR_b$ is the fiber
above $b$ of the path connected component of $S_\sigma$.

For each $\sigma \in \widehat G$, we define then $H_\sigma^i$ to be the
union of the isotypical components corresponding to $S_{\sigma,b}$ in
$H_b^i$. Then $H_\sigma^i$ is again a bundle of Hilbert
spaces. Moreover, $F_b$ will map $H_\sigma^0$ to $H_\sigma^1$, and the
resulting map will be an isomorphism except maybe for finitely many
irreducible representations $\sigma$.

We shall construct the bundle $\KER$ as a union of sub-bundles
corresponding to each representation $\sigma \in \widehat G$, for those
$\sigma \in \widehat G$ for which $F_b : H_\sigma^0 \to H_\sigma^1$ is not
already an isomorphism. We can thus fix $\sigma$ in the following
discussion.

If the bundle $H_\sigma^0$ is consists of finite dimensional vector
spaces, then $H_\sigma^1$ consists also of finite dimensional vector
spaces, and hence we can choose $\KER = H_\sigma^0$.  Let us assume
then that $H_\sigma^0$ does not consist of finite dimensional vector
spaces, then $H_\sigma^1$ does not consist either of finite
dimensional vector spaces. The triviality of the bundle $H_\sigma^0$
implies then that we can choose an increasing sequence of finite
dimensional sub-bundles $L_n \subset H_\sigma^0$ such that their union
is dense in each fiber of $H_\sigma^0$. Using the fiberwise averaging
with respect to $\GR_b$, we see that we can assume each $L_n$ to be
invariant. We can then take the component of $\KER$ corresponding to
$\sigma$ to be $L_n$, for some large $n$.
\end{proof}\smallskip

\begin{dfn}\label{def.giindex}\
Let $F=(F_b :H^0_b \to H^1_b)_{b \in
B}$ be a family of $\GR$-invariant Fredholm operators that is norm-continuous in
any trivialization of $H^i \to B$, as in Lemma \ref{lemma.KER}. The element
$$
        \indg (F):=[\KER]-[\COK]\in K_\GR(B)
$$
is called the {\em gauge-equivariant index\/} of the invariant
Fredholm family $F$.
\end{dfn}


\begin{lem}\label{lem:indexindep}
The gauge-equivariant index $\indg F$ is well defined, that is, it
depends only on $F$ and not on the choice of the $\GR$-equivariant
sub-bundle $\KER$.
\end{lem}

\begin{proof}\
Let us consider two possible choices for the bundle $\KER$ of Lemma
\ref{lemma.KER} that was used to define the gauge-equivariant index.
Denote these two sub-bundles by $\KER_1$ and $\KER_2$ and identify
them with the orthogonal the projections onto their ranges. Let $P_1$
and $P_2$ be these two projection. The proof of Lemma \ref{lemma.KER}
shows that we can find a new sub-bundle $\KER$ with associated
projection $P$ such that $\|P_1 - P_1P\| < \epsilon$ and $\|P_2 -
P_2P\| < \epsilon$, with $\epsilon$ as small as we want. It is enough
to check then that both $\KER_1$ and $\KER$ give rise to the same
index.

But $P_1$ is close to the subprojection $\chi(PP_1P)$ of $P$ obtained
by applying the functional analytic calculus to $PP_1P$, where $\chi$
is locally constant, equal to $0$ in a neighborhood of $0$ and equal
to $1$ in a neighborhood of $1$. Since close projections are homotopic
and the index does not change under homotopies, we see that we can
actually assume that $P_1 \le P$ (that is, that $P_1$ is a
subprojection of $P$). But then $F$ is injective from the range of $P
- P_1$, $(P-P_1)H^0$ to $F(P-P_1)H^1$. Moreover, $F(P - P_1)H$ is also
a finite-dimensional vector bundle over $B$. Let $\COK_1$ be the
cokernel of $F$ acting on $P_1H^0$. The above discussion shows that
$$
        [\KER] - [\COK] = [{\KER}_1] + [(P-P_1)H^0] - [{\COK}_1] -
        [F(P-P_1)H^1] = [{\KER}_1] - [{\COK}_1].
$$
\end{proof}\smallskip

The gauge-equivariant index has the usual properties of the index of
elliptic operators.

\begin{lem}\label{lem:indexbliz}
Let $F$ and $F'$ be two invariant families as in Lemma
\ref{lemma.KER}. Assume $F$ consists of Fredholm operators. Then there
exists $\epsilon > 0$ such that if $\|F - F'\| < \epsilon$, then $F'$
consists also of Fredholm operators and has the same gauge-equivariant
index. In particular, the gauge-equivariant index is homotopy
invariant.
\end{lem}

\begin{proof}\
The family $F'$ is Fredholm by the usual Hilbert space argument which
applies since our operators are norm continuous in any
trivialization. Moreover, for $F'$ sufficiently close to $F$, we can
choose the same sub-bundle $\KER$ to define the gauge-equivariant
index of $F'$, while the corresponding $\COK$ and $\COK'$ will be
isomorphic. See also \cite{TroVINITI99}.
\end{proof}\smallskip

Let us now consider a \LS $\GR$-fiber bundle $Y \to B$ with a
$\GR$-invariant complete metric on each fiber of $Y$. The invariant
metrics give rise to Laplace-Beltrami operators $\Delta_b$ acting on
functions.  The Sobolev spaces $H^l(Y)$ are defined as the domains of
$\Delta_b^l$ (this choice is classical, see \cite{Aubin, Roe}, for
example).  This definition extends right away to Sobolev spaces of
sections of a $\GR$-equivariant vector bundle $E$ equipped with a
$\GR$-invariant hermitian metric to define a locally trivial bundle of
$\GR$-Hilbert spaces.

Assume that $F = (D_b)_{b \in B}$ is a $\GR$-invariant family of order
$m$ operators acting between sections of some $\GR$-equivariant vector
bundles $E_0,E_1 \to Y$. The $l$th Sobolev spaces $H^l(Y;E_0)$ of
sections of $E_0$ along the fibers of $Y \to B$ defines a locally
trivial bundle of $\GR$-Hilbert spaces $H_0$. Let $H_1$ be defined
similarly using the $(l-m)$th Sobolev spaces of sections of
$E_1$. Then $F : H_0 \to H_1$ is a Fredholm family as in the statement
of Lemma \ref{lemma.KER}. The homotopy invariance of the
gauge-equivariant index shows that the gauge-equivariant index of the
family $F = (D_b)_{b \in B}$ depends only on the principal symbol of
this family. As in the non-equivariant case, the principal symbol of
$F$ gives rise to an element $x$ of $K^0_{\GR}(T^*Y_{vert})$, and the
index depends only on $x$. Since every class $x \in
K^0_{\GR}(T^*Y_{vert})$ arises in this way, we obtain a well defined
group morphism
\begin{equation}\label{eq.a.i.g}
        a-\ind: K^0_{\GR}(T^*Y_{vert})\to K^0_\GR(B),
\end{equation}
which will be called the {\em analytic index morphism}.  A more
general definition of this morphism (without finite holonomy
conditions will be obtained in Section \ref{S.An.Ind}.

The locally trivial bundle of $\GR$-Hilbert spaces $H_0$ defined above
will be called the bundle of Sobolev spaces of order $l$ associated to
$Y \to B$.

\begin{dfn}
A locally trivial bundle of $\GR$-Hilbert spaces $\pi:H \to B$ is
called {\em saturated} if, for any $b \in B$ and any $\sigma \in
\widehat {\GR_b}$, the multiplicity of $\sigma$ in the Hilbert space
$H_b=\pi^{-1}(b)$ is either zero of infinite.
\end{dfn}

One has the following easy but useful statement.

\begin{lem}\label{lem:dimsatur}
Suppose that $\dim Y > \dim\GR$. Than any bundle of Sobolev spaces
$H^s(Y;E)$ associated to $Y \to B$ is saturated.
\end{lem}

\begin{proof}\
Let $G$ be the typical fiber of $\GR \to B$ and fix $\sigma \in
\widehat G$.  Let $H \to B$ be a bundle of Sobolev spaces associated
to $Y \to B$.  Then the multiplicity of $\sigma$ in the fiber $H_b$ is
a multiple of the dimension of $L^2(Y_b/\GR_b)$.
\end{proof}\smallskip

The following lemma explains why we are interested in saturated
Hilbert bundles.

\begin{lem}\label{lem:saturvlozh}
Suppose that $H \to B$ is a saturated locally trivial bundle of
$\GR$-Hilbert spaces. Let $E \to B$ be a $\GR$-equivariant vector
bundle. Assume that for any $b \in B$ and any $\sigma \in \widehat
{\GR_b}$ appearing in $E_b$, the multiplicity of $\sigma$ in $H_b$ is
non-zero. Then $E$ is $\GR$-equivariantly isomorphic to a sub-bundle
of $H$.
\end{lem}

\begin{proof}\
Choose compact subsets $U_\alpha \subset B$ and trivializations
$$
        \GR\vert_{U_\alpha} \simeq U_\alpha \times G, \quad
        E\vert_{U_\alpha} \simeq U_\alpha \times E_0, \quad \text{and
        }\; H\vert_{U_\alpha} \simeq U_\alpha \times H_0.
$$
We can choose the sets $U_\alpha$ such that their interiors cover $B$.

We can choose then embeddings $J_\alpha : E\vert_{U_\alpha} \to
H\vert_{U_\alpha}$ inductively such that $J_\alpha$ and $J_\beta$ have
orthogonal ranges above each $b \in U_\alpha \cap U_\beta$, if $\alpha
\not = \beta$.  Let $\phi_\alpha^2$ be a partition of unity
subordinated to the interiors of $U_\alpha$. Then $J = \sum
\phi_\alpha J_\alpha$ is the desired embedding.
\end{proof}\smallskip

For the proof of the next theorem we shall need the following lemma.

\begin{lem}\label{lem:stabildliasatur}
Let $E_1, E_2 \to B$ be two $\GR$-equivariant vector bundles on a
path-connected, locally simply-connected topological space $B$ such
that their classes in gauge equivariant theory coincide (that is,
$[E_1] = [E_2] \in K^0_\GR(B)$).  Then there exists a $\GR$-fiber
bundle $E \to B$ such that
\begin{enumerate}
\item $E_1 \oplus E \cong E_2 \oplus E$;
\item if an irreducible representation $\sigma \in \widehat{\GR}_b$
appears in $E_b$, then it appears also in $(E_1)_b$ and in $(E_2)_b$.
\end{enumerate}
\end{lem}

\begin{proof}\
The existence of $E$ satisfying the first assumption follows from the
definition of the group completion of a monoid.

To obtain the second property, we just decompose the vector bundles
according to representations $\sigma$ in the orbits of $\pi_1(B,b_0)$
on $\widehat G$. (See the discussion at the end of Section
\ref{S.Inv.Ops}.) Then we notice that the dimension of the isotypical
subspace $(E_b)_\sigma$ is the same for any $\sigma \in
\widehat{\GR}_b$ belonging to a fixed connected component of
$\widehat{\GR}$.
\end{proof}\smallskip

In the proof of the following theorem, we shall use the completion of
the algebra $\Psm{-\infty}{Y;E}$ of order $-\infty$, invariant
operators on $Y$ and acting on square integrable sections of $E$. We
shall denote by $C^*(Y;\GR,E)$ the resulting algebra. It has a natural
norm
\begin{equation}\label{eq.def.norm}
        \| T \| = \sup_{b \in B} \|T_b\|_b,
\end{equation}
where $\|\;\|_b$ is the norm of operators acting on the Hilbert space
$L^2(Y_b;E_b)$, $E_b := E\vert_{Y_b}$. Since all operators in
$\Psm{-\infty}{Y;E}$ act as compact operators on $L^2(Y_b;E_b)$, a
density and an averaging argument shows that $C^*(Y;\GR,E)$ identifies
with the algebra of continuous families of $\GR$-invariant, compact
operators acting on the family of Hilbert spaces $L^*(Y_b;E_b)$. If
$E$ is a trivial vector bundle, we omit it from the notation.  Also,
if $Y = \GR$, then we shall denote the resulting algebra $C^*(Y;\GR)$
simply by $C^*(\GR)$. Recall for the next theorem that the spaces
$\Psm{m}{Y;E,F}$ were defined in Section \ref{S.Inv.Ops}.

The following theorem shows that the $\GR$-equivariant index identifies
the obstruction to invertibility, as the usual (or Fredholm) index.

\begin{teo}\label{thm.pert.I}\
Suppose that $\dim Y>\dim\GR$ and let $D \in \Psm {m}{Y;E,F}$ be a
$\GR$-equivariant family of elliptic operators acting along the fibers
of $Y \to B$.  Then we can find $R \in \Psm {m-1}{Y;E,F}$ such that
$$
        D_b + R_b : H^{s}(Y_b;E_b) \to H^{s-m}(Y_b;F_b)
$$
is invertible for all $b \in B$ if, and only if, $\indg(D)=0$.
Moreover, if $\indg(D)=0$, then we can choose the above $R$ in
$\Psm{-\infty}{Y;E,F}$.
\end{teo}

\begin{proof}\
Suppose that such a perturbation $R$ exists. Then $\indg (D + R) = 0$.
Since the index depends only on the principal symbol of $D$, we obtain
that $\indg(D)=0$.

Let $D \in \Psm{m}{Y;E,F}$ and choose $\KER$ and $\COK$ as in Lemma
\ref{lemma.KER}. If $\indg(D)=0$, then, by Lemma
\ref{lem:stabildliasatur}, we can find a $\GR$-equivariant vector
bundle $E \to B$ such that $\KER \oplus E \cong \COK \oplus E$. We can
choose this bundle $E$ such that all irreducible representations
$\sigma \in \widehat{\GR}$ that appear in (some fiber of) $E$ also
appear in (some fiber of) the bundle of Sobolev spaces $H^s(Y;E)$.
Lemma \ref{lem:saturvlozh} then shows that we can identify $E$ with a
$\GR$-equivariant sub-bundle the orthogonal complement of $\KER$.
Then, by replacing our original choice of $\KER$ with $\KER \oplus E$
and then by replacing $\COK$ with the new cokernel space, we can
assume that $\KER \cong \COK$.

Let $T : \KER \to \COK$ be a $\GR$-equivariant isomorphism of these
two $\GR$-equivariant vector bundles.  Then $T$ is a $\GR$-invariant
family of compact operators acting on sections of $H^s(Y;E)$ with
values sections of $H^{s-m}(Y;F)$. Let $P$ be the orthogonal
projection onto $\KER$. Then $D'=D(1-P) + T$ is invertible. Moreover,
$R' := D' - D$ is also a $\GR$-invariant family of compact operators
acting on sections of $H^s(Y;E)$ with values sections of
$H^{s-m}(Y;F)$. Since $\Psm{-\infty}{Y;E}$ is dense in $C^*(Y;\GR,E)$,
we can find an operator $R \in \Psm{-\infty}{Y;E}$ that is close
enough to $R'$ to ensure that $D + R$ is also invertible in all
fibers.
\end{proof}\smallskip

 \section{The structure of $\GR$-equivariant $K$-groups\label{Sec.structure.K}}

In this section we shall study the structure of the algebra
$C^*(\GR)$ introduced in the previous section as the
completion of $\Psm{-\infty}{\GR}$ in the norm $\| \; \|$ of Equation
\ref{eq.def.norm}. We shall also relate the gauge-equivariant
$K$-theory groups of $\GR$ with the $K$-theory groups of $C^*(\GR)$.

For each $b \in B$, denote by $\mu_b$ the translation invariant
measure on $\GR_b$ whose total mass is one. Because the bundle $\GR
\to B$ is locally trivial, we know that the function $B \ni b \to
\mu_b(f) \in \CC$ is continuous, for any continuous function $f$ on
$B$. Let us denote by $C(B)$ the space of continuous functions on
$\GR$ with the fiberwise convolution product and the involution
$f^*(g) = \overline{f(g^{-1})}$. The algebra $C(\GR)$ also acts on
each $L^2(\GR)$ and we can check using the local trivialization of
$\GR$ that $C(\GR)$ is dense in $C^*(\GR)$.

The algebra that we have introduced above is usually denoted
$C^*_r(\GR)$, whereas the notation $C^*(\GR)$ is reserved for the
envelopping $C^*$-algebra of $C(\GR)$. It can be shown, but we shall
not need this, that in our case $C^*(\GR) = C^*_r(\GR)$, which
justifies our notation. (See \cite{LN} for the definition of the
envelopping $C^*$-algebra of a groupoid and for related concepts).

Recall that we have denoted by $\mP$ a principal $\Aut(G)$-bundle that
defines $\GR$ in the sense that $\GR = \mP \times_{\Aut(G)} G$, for
some fixed Lie group $G$. Also, recall that in this paper $G$ was
assumed to be compact beginning with Section \ref{Sec.Gauge.equiv}.

If $\GR$ (or, equivalently, $\mP$) is trivial, then $C^*(\GR) \simeq
C(B, C^*(G))$, the algebra of continuous functions on $B$ with values
in $C^*(G)$. (As usual, we have denoted by $C^*(G)$ the norm completion
of the convolution algebra of $G$
actiong on $L^2(G)$.)  This leads us to the following construction.
Let $C^*_\GR\to B$ be the locally trivial bundle with fibers
$C^*(\GR_b)$.  It is the fiber bundle associated to $\mP$ and its
action on $C^*(G)$, that is, $C^*_\GR \cong \mP \times_{\Aut(G)}
C^*(G)$ as bundles over $B$.
The local triviality of this bundle allows us to talk about
the space of continuous sections of this bundle, which is a complete
normed algebra (even a $C^*$-algebra).

\begin{lem}
The algebra $C^*(\GR)$ identifies naturally with the algebra
$\G(C^*_\GR)$ of continuous sections of the bundle $C^*_\GR \to B$.

If $\GR$ has representation theoretic finite holonomy, then there
exist projections $p_n \in C^*(\GR)$, such that $p_n p_{n+1} = p_n$,
$p_n x = x p_n$, for any $x \in C^*(\GR)$, and $\cup p_n C^*(\GR)$ is
dense in $C^*(\GR)$.
\end{lem}

\begin{proof}\
The first statement follows from the fact that these two algebras are
the completions of the same algebra $\Psm{-\infty}{\GR}$ with respect
to the same norm.

Assume first that $B$ is path-connected. We shall use the notation and
the constructions introduced at the end of Section \ref{S.Inv.Ops}.
In particular, $H_\GR$ is the image of the holonomy morphism
$\pi_1(B,b_0) \to H_{R} := \Aut(G)/\Aut_0(G)$. Let $H$ be the inverse
image of $H_{\GR}$ in $\Aut(G)$. Then we can reduce the structure
group of $\mP$ to $H$.  Choose sets $S_n \subset \wh G$ such that $S_n
\subset S_{n+1}$, each $S_n$ is $\pi_1(B,b_0)$-invariant and $\cup S_n
= \widehat G$.  Also, let $q_n$ be the central projection of $C^*(G)$
corresponding to $S_n$.  By construction, $q_n$ is invariant for $H$,
and hence it gives rise to a section $p_n$ of $C^*_\GR$, which is the
desired projection.

In general, we use an exhaustion of $\widehat{\GR}$ by compact-open
subsets.
\end{proof}\smallskip

Let $A$ be a (possibly non-unital) algebra. By a {\em
finitely-generated, projective module} over $A$ we shall understand a
left $A$-module of the form $A^Ne$, where $e \in M_N(A)$ is a
projection (that is, $e^2 = e$).  All modules over non-commutative
algebras used below will be assumed to be left-modules, unless
otherwise mentioned.

\begin{teo}\label{thm.equiv.cat}\ Assume that the bundle of
compact Lie groups $\GR \to B$ has representation theoretic finite
holonomy and that $B$ is compact.  Then there is a natural equivalence
of categories between the category of locally trivial
$\GR$-equivariant vector bundles over $B$ and the category of
finitely-generated, projective modules over $C^*(\GR)$.  In
particular,
$$
        K^*_\GR(B) \cong K_*(C^*(\GR)).
$$
\end{teo}

\begin{proof}\
If $E \to B$ is a $\GR$-equivariant vector bundle, we can endow $E$
with a $\GR$-invariant metric (Corollary \ref{cor.herm}) and hence we
obtain that $C(\GR)$ acts on $\Gamma(E)$. Using the local triviality
of $\GR$ and the metric on $E$, we see that we can extend this action
of $C(\GR)$ to an action of $C^*(\GR)$. Thus $\Gamma(E)$ is a
$C^*(\GR)$-module. We shall prove that if is projective and that
$E \to \Gamma(E)$ is an equivalence of categories.

By looking at the representations of $\GR$ that appear in the fibers
of $E \to B$, we see that by choosing $n$ large enough we can assume
that $p_n$ acts as the identity on $\Gamma(E)$. Then there exists a
surjective map $(C^*(\GR)p_n)^N \to \Gamma(E)$ of
$C^*(\GR)$-modules. By regarding $\cup_b (C^*(\GR_b)p_n)^N$ as a
$\GR$-equivariant vector bundle, we see that $E$ has a direct summand
in it, which implies then that the $C^*(\GR)$-module $\Gamma(E)$ is a
direct summand of $(C^*(\GR)p_n)^N$, and hence it is of the form
$C^*(\GR)^Ne$, for some projection $e \in M_N(p_nC^*(\GR)p_n)$.

Conversely, let us assume that $M$ is a finitely-generated, projective
$C^*(\GR)$-module, that is, that $M \cong C^*(\GR)^N e$, for some
projection $e \in M_N(C^*(\GR))$. Then the space of continuous sections
of $\GR$ will also act on $M$. Using the local triviality of $\GR \to
B$, we then see that there exists $n$ such that $p_n e = e$. Then $M =
(p_n C^*(\GR))^N e$ is also a projective module over the unital algebra
$p_n C^*(\GR) p_n = p_n C^*(\GR)$, which contains $C(B)$ in its
center. Since $p_n C^*(\GR)$ is a projective $C(B)$-module, $M$ will
also be a projective $C(B)$-module, and hence it can be identified
with the space of sections of a vector bundle $E \to B$.  Then the
action of $C^*(\GR)$ on $E$ is fiberwise and, by restricting to $\GR$,
we obtain that $E$ is also a $\GR$-equivariant vector bundle.
Thus $M \cong \Gamma(E)$ as $C^*(\GR)$-modules.
\end{proof}\smallskip

We remark that the first part of the above theorem remains true even
without the representation theoretic finite holonomy condition, but
the proof has to be slightly modified. The second part, that is the
isomorphism of the $K$-theory groups, is not true, in general, without
the representation theoretic finite holonomy condition.

We now take a closer look at the structure of the algebra $C^*(\GR)$.
Let us denote by $(\widehat{\GR})_d$ the space of irreducible
representations of dimension $d$ of the groups $\GR_b$. By the local
triviality of $\GR \to B$, $(\widehat{\GR})_d$ is open and closed in
$\widehat{\GR}$ and is a covering space of $B$. Let $f$ be a
continuous function on $\GR$. Then the function
$$
        (\widehat{\GR})_d \ni \sigma \to f_{TR} := Tr(\sigma(f)) \in
        \CC
$$
is a continuous function on $(\widehat{\GR})_d$ (this also follows
from the local triviality of $\GR \to B$). Moreover,
$$
        |Tr (\sigma (f - g))| \le (\dim \sigma)\, \|f - g \|,
$$
so if $f_n \in C(\GR)$ converges in $C^*(\GR)$, the functions
$(f_n)_{TR}$ converge uniformly on each of the sets
$(\widehat{\GR})_d$, which shows that the definition of $f_{TR}$ can
be extended to $f \in C^*(\GR)$ by continuity, and the result will
still be continuous on $\widehat{\GR}$.

Let us denote by $PGL(d,\CC) := GL(d,\CC)/Z(GL(d,\CC))$ the group of
automorphisms of the algebra $M_d(\CC)$. If $\mathcal A \to X$ is a
bundle of algebras with structure group $PGL(d,\CC)$, then every fiber
$\mathcal A_x \cong M_d(\C)$ will have a unique $C^*$-norm denoted $\|
\; \|_x$, for any $x \in X$.  (Recall that a norm $\|\; \|$ on a
$*$-algebra $A$ is called a $C^*$-norm if it is a complete Banach
algebra norm and $\|x^*x\| = \|x\|^2$ for any $x \in A$. A normed
$*$-algebra is called a {\em $C^*$-algebra} if its norm is a
$C^*$-norm.)  We shall denote by $\Gamma_0(\mathcal A)$ the space of
continuous sections $\xi$ of $\mathcal A$ such that for any $\epsilon
> 0$, the set $\{ x, \| \xi \|_x \ge \epsilon \}$ is a compact subset
of $\widehat{\GR}$. Then $\Gamma_0(\mathcal A)$ is complete in the
norm $\|\xi \| = \sup_x \|\xi(\sigma)\|_x$ and is a $C^*$-algebra.

\begin{teo}\ Let $\GR \to B$ be a bundle of compact Lie
groups. Then there exists on each $(\widehat{\GR})_d$ a locally
trivial bundle of algebras $\mathcal A_d$ with fiber $M_d(\CC)$ and
structure group $PGL(d,\CC) := GL(d,\CC)/Z(GL(d,\CC))$ such that the
space $\Gamma_0(\mathcal A_d)$ identifies with a direct summand of
$C^*(\GR)$ and
$$
        C^*(\GR) \cong \Gamma_0(\mathcal A_d).
$$
In particular, $K_i(C^*(\GR)) \cong \oplus K_i(\Gamma_0(\mathcal
A_d))$ and the primitive ideal spectrum of $C^*(\GR)$ is homeomorphic
to $\widehat{\GR}$, which in turn is homeomorphic to the disjoint
union of the sets $(\widehat{\GR})_d$.
\end{teo}

\begin{proof}\
The first part of the result follows from the fact that the pointwise
trace $f \to f_{TR}$ defined above is continuous for any $f \in
C^*(\GR)$. The second part is also a general property of continuous
trace $C^*$-algebras. (See \cite{Dixmier} and the references
therein.)
\end{proof}\smallskip

We thus obtain a description of $K_i(C^*(\GR))$ in terms of twisted
$K$-theory. (A {\em twisted $K$-theory group} of a space $X$ is the
$K$-theory of a bundle of matrix algebras or compact operators over
$X$. See \cite{Mathai, Rosenberg, FW} and the references therein for
more information on the subject.)

The above theorem has several consequences.

\begin{prop}\
Suppose that all the fibers of $\GR \to B$ are compact abelian Lie
groups. Then
$$
        K_*(C^*(\GR)) \cong \oplus K^*(\widehat{\GR}).
$$
\end{prop}

\begin{proof}\
If the fibers of $\GR \to B$ are abelian groups, then $C^*(\GR)$ is
also abelian and its primitive ideal spectrum is $\widehat{\GR}$.
\end{proof}\smallskip

We shall assume until the end of this section
that $B$ is a path-connected, locally simply-connected space.
(Also, recall that beginning with Section \ref{Sec.Gauge.equiv},
$B$ is assumed to be compact.)

Let $\mP$ be the principal $\Aut(G)$ bundle defining $\GR \to B$.  We
shall identify $G$ with one of the fibers of $\GR \to B$.  Fix $\sigma
\in \widehat G$, we shall denote then by $B_\sigma$ the connected
component of $\widehat{\GR}$ containing $\sigma$.  Recall that
$\pi_1(B,b_0)$ acts on $\widehat G$. Then $B_\sigma \to B$ is the covering
space associated to the isotropy of $\sigma$ in $\pi_1(B,b_0)$, more
precisely, if $\tilde B \to B$ is a universal covering space of $B$,
then the covering $B_\sigma \to B$ is equivalent to the covering
$\tilde B \times_{\pi_1(B,b_0)} (\pi_1(B,b_0)\sigma) \to B$.

The space $B_\sigma$ is contained in the space $(\widehat{\GR})_d$, with
$d = \dim \sigma$. We shall denote by $\mathcal A_\sigma$ the bundle
of finite dimensional algebras obtained by restricting $\mathcal A_d$
to $B_\sigma$. We then have the following corollary.

\begin{cor}\label{cor.direct.sigma}\
Assume $\GR \to B$ is a bundle of compact Lie groups over a
path-connected, locally simply-connected space.
Then we have a canonical isomorphism
$$
        C^*(\GR) \cong \oplus_{\sigma} \Gamma_0(\mathcal A_\sigma),
$$
where $\sigma$ ranges through a set of representatives of the orbits
of $\pi_1(B,b_0)$ on $\widehat G$ and $\mathcal A_\sigma \to B_\sigma$ are
bundles of algebras with fiber $M_d(\CC)$, $d = \dim \sigma$, obtained
as the restriction of $\mathcal A$ to $B_\sigma$, as above.
\end{cor}

For the following result, we shall need the following lemma.

\begin{lem}\label{lemma.rat.iso}\ Let $\mathcal A \to X$ be
a locally trivial bundle of finite dimensional algebras with fiber
$M_d(\CC)$. Then $K_i(\Gamma_0(\mathcal A)) \otimes \QQ \cong K^i(X)
\otimes \QQ$.
\end{lem}

\begin{proof}\
Let $A = \Gamma_0(\mathcal A)$. Endow each fiber of $\mathcal A \to X$
with the inner product $(T,T_1) := Tr(T_1^*T)$, with $T$ and $T_1$ in
the fiber. Then $\mathcal A \to X$ is a Hermitian vector bundle. Let
$B = \Gamma_0(\End(\mathcal A))$, which will be then a $C^*$-algebra
Morita equivalent to $C_0(X)$.  In particular, $K_i(C_0(X)) \cong
K_i(\Gamma_0(\End(\mathcal A)))$, isomorphism that we shall denote by
$j$ in what follows. Note also that we have a natural morphism $j_0 :
C_0(X) \to B$, which sends a function $f$ to the operator of
multiplication with that function $f$. Let $[\mathcal A]$ be the class
in $K^0(X) = K_0(C_0(X))$ of the vector bundle $\mathcal A \to
X$. Then, at the level of $K_0$-groups, we have $(j_0)_*(\xi) =
j_*(\xi) \otimes [\mathcal A]$.

The center of $A$ is isomorphic to $C_0(X)$, which gives rise to an
algebra morphism $\phi : C_0(X) \to A$. Left multiplication with the
elements of $A$ gives rise to a second algebra morphism $\psi : A \to
B$.  The composite morphisms $(j \circ \psi_*) \circ \phi_* :
K_0(C_0(X)) \to K_0(B) \simeq K_0(C_0(X))$ is then multiplication by
$[\mathcal A]$.

Assume $\mathcal A$ is a trivial bundle of rank $d^2$. Then $(j \circ
\psi_*) \circ \phi_* : K_0(C_0(X)) \to K_0(B) \simeq K_0(C_0(X))$ is
multiplication by $d^2$.  Similarly, $\phi_* \circ (j \circ \psi_*) :
K_0(A) \to K_0(A)$ is also multiplication by $d^2$.  This proves the
rational isomorphism for trivial bundles $\mathcal A$.

In general, we use the fact that the maps $\phi_*$ and $\psi_*$ are
natural and a Meyer-Vietoris argument \cite{Blackadar}, to complete
the proof.
\end{proof}\smallskip

The above example is not true if we do not include rational
coefficients, as implied by the following example due to M. Dadarlat.

\newtheorem{example}[teo]{Example}

\begin{example}{\rm
Let $X$ be the mapping cone of the map $z \mapsto z^n$ of the
circle. Then $X$ is a two dimensional CW complex with
\[K_0(C(X))=\mathbb{Z} \oplus \mathbb{Z}/n, \,\,K_1(C(X))=0.\]

Using the homotopy exact sequence
\[\mathbb{T} \rightarrow U(n)\rightarrow
PU(n)\rightarrow B\mathbb{T}\rightarrow BU(n),\] we see that any
$C(X)$-linear automorphism of \[A=C(X)\otimes M_n(\mathbb{C})\] that
is given by a map $X \rightarrow PU(n)$ is determined up to unitary
equivalence by a line bundle $E$ over $X$ with $E \oplus \cdots \oplus
E$ ($n$-times) trivial.

One also can verify that the map $\alpha_*:K_0(A) =K^0(X)\rightarrow
K_0(A)=K^0(X)$ is induced by multiplication (in the ring $K^0(X))$
with the $K_0$-class $[E]$.  Define the mapping cylinder $M_\alpha$ of
$\alpha$:
\[M_\alpha=\{f:[0,1] \to A: f(1)=\alpha(f(0))\}.\]
Then $M_\alpha$ is the $C^*$-algebra of sections of a bundle of $n
\times n$-matrices with spectrum $Y=\mathbb{T}\times X$ and center
$C(Y)$. The exact sequence
\[0 \rightarrow SA \rightarrow M_\alpha \rightarrow A \rightarrow 0\]
induces a six-term exact sequence
$$
	\CD 0=K_1(A) @>>>K_0(M_\alpha) @>>>K_0(A)=\mathbb{Z} \oplus
	\mathbb{Z}/n\\ @A{1-\alpha_*}AA @. @VV{1-\alpha_*}V\\ 0=K_1(A)
	@<<< K_1(M_\alpha) @<<< K_0(A)=\mathbb{Z} \oplus \mathbb{Z}/n
	\endCD
$$
Choose $E$ such that its class in $K^0(X)\cong \mathbb{Z} \oplus
\mathbb{Z}/n$ is equal to $(1,-1)$ and set $n=4$. Then $K_0(M_\alpha)=
\mathbb{Z} \oplus \mathbb{Z}/2$ so that $K_0(M_\alpha)\neq
K_0(C(Y))=K_0(C(\mathbb{T}\times X))=\mathbb{Z} \oplus
\mathbb{Z}/4$. This completes the example of a bundle of algebras
whose $K$-theory is not isomorphic to that of its center.}
\end{example}

The decomposition of Corollary \ref{cor.direct.sigma}
leads to the following determination of
the $K$-theory groups of the algebras $C^*(\GR)$, up to
rational isomorphism.

\begin{teo}\
Suppose $B$ is a path-connected, locally simply-connected space and
let $\GR \to B$ be a bundle of compact Lie groups with representation
theoretic finite holonomy. Then
$$
        K_i(C^*(\GR)) \otimes \QQ \simeq K^i(\widehat{\GR}) \otimes
        \QQ \simeq K^i(B) \otimes R(G)^{\pi_1(B,b_0)} \otimes \QQ.
$$
\end{teo}

\begin{proof}\
Let $G$ be the typical fiber of $\GR \to B$ and choose a set of
representatives $S \subset \wh G$ for the orbits of $\pi_1(B,b_0)$ on
$\wh G$. We shall use the results and the notation of Corollary
\ref{cor.direct.sigma}.

We have then
$$
        K_i(C^*(\GR)) \otimes \QQ \simeq \oplus_{\sigma \in S}
        K_i(\Gamma_0(\mathcal A_\sigma)) \otimes \QQ \simeq
        \oplus_{\sigma \in S} K^i(B_\sigma) \otimes \QQ.
$$
But $K^i(B_\sigma) \otimes \QQ \simeq K^i(B) \otimes \QQ$, because
$B_\sigma \to B$ is a finite covering. Thus $K_i(C^*(\GR)) \otimes \QQ
\simeq K^i(B) \otimes \QQ^{S}$, where $\QQ^{(S)}$ is the vector space
with basis $S$. Moreover, $\QQ^{(S)}$ identifies with
$R(G)^{\pi_1(B,b_0)} \otimes \QQ$, because $R(G) = \ZZ^{(\wh G)}$
(that is, the free abelian group with basis $\wh G$).
\end{proof}\smallskip

Recall that the conditions of the above theorem are automatically
satisfied if the typical fiber $G$ of $\GR \to B$ is a semi-simple Lie
group or if $G$ is the product of a semisimple Lie group by the
one-dimensional torus $S^1$.

Let us now take a closer look at the algebras $A_\sigma =
\Gamma_0(\mathcal A_\sigma)$ used above. Each of these algebras is the
algebra of sections of the field $\mathcal A_\sigma$ of finite
dimensional matrix algebras over $B_\sigma$. Let us denote the
dimension of these fibers by $d_\sigma^2$, for any fixed $\sigma$, as
above. In particular, $d_\sigma = \dim V_\sigma$. Then the bundle
$\mathcal A_\sigma$ is a bundle with structure group
\begin{equation}
        PGL(d_\sigma,\CC) : = GL(d_\sigma,\CC)/ Z(GL(d_\sigma,\CC)) =
        SL(d_\sigma,\CC)/ C_{d_\sigma} =: PSL(d_\sigma,\CC),
\end{equation}
where by $C_m$ we denote the cyclic group with $m$ elements.  The
bundle $\mathcal A_\sigma$ is hence classified by a 1-cocycle in
$H^1(B_\sigma, PSL(d_\sigma, \CC))$. The connecting morphism
$$
        H^1(B_\sigma, PSL(d_\sigma, \CC)) \to H^2(B_\sigma,C_{d_\pi})
$$
in non-abelian cohomology then gives rise to an element $\chi_\sigma
\in H^2(B_\sigma,C_{d_\pi})$ (see \cite{Dixmier, DonovanKaroubi}).
By definition, $\chi_\sigma$ is the Dixmier-Duady invariant of
$\mathcal A_\sigma$ (see \cite{Dixmier, DonovanKaroubi}). We want to
analyse these invariants for the fields of matrix algebras $\mathcal
A_\sigma \to B_\sigma$ introduced above before the statement of
Corollary \ref{cor.direct.sigma}.

Assume now that $G$, the typical fiber of $\GR$, is connected.  We
shall use the notation used at the end of Section \ref{S.Inv.Ops}.  In
particular, $\mP$ is the principal $\Aut(G)$-bundle defining $\GR$ and
$$
        \pi_1(B,b_0) \to H_{R} = \Aut(G)/\Aut_0(G)
$$
is the holonomy morphism defining the principal $H_R$-bundle $\mP_0 :=
\mP/\Aut_0(G)$.  Recall that, in our case ($G$ connected), $\GR \to B$
has representation theoretic finite holonomy if, and only if, the
connected components of $\widehat \GR$ are compact.  We then know by
Theorem \ref{thm.G=conn} that these conditions are also equivalent to
the condition that the range $H_\GR$ of the holonomy morphism be
finite. We identify $G$ with a fiber of $\GR \to B$ and let $B_\sigma$
be the connected component in $\widehat \GR$ of $\sigma \in \widehat
G$, as above.  Let $\tilde B \to B$ be a universal covering space of
$B$ and let $B' := \tilde B \times_{\pi_1(B,b_0)} H_R$. Then for any
$\sigma$ we obtain a map
\begin{equation}\label{eq.A}
        f_\sigma : B' \to B_\sigma = \tilde B \times_{\pi_1(B,b_0)}
        H_R\sigma.
\end{equation}
Let us lift $\mathcal P$ to an $\Aut(G)$-principal bundle on
$B'$. Then, by construction, the resulting bundle reduces to a
principal $\Aut_0(G) = G_{int} = G/Z(G)$ bundle classified by a
one-cocycle in $H^1(B', G_{int})$ (we are using here Corollary
\ref{prop.Z=1}). Let us denote by $\chi \in H^2(B', Z(G))$ the image
of this cocycle under the connecting morphism in non-abelian
cohomology associated to the exact sequence of groups
$$
	1 \to Z(G) \to G \to G/Z(G) \to 1.
$$  
Let $G'$ be the derived group of $G$ ($G'$ is connected because we
assumed $G$ to be connected).  Denote
\begin{equation}\label{eq.B}
	Z' = G' \cap Z(G),
\end{equation}
a finite abelian group.  Because $G'$ maps onto $G/Z(G)$, the
obstruction $\chi$ comes from a canonical element
\begin{equation}\label{eq.chi}
        \chi' \in H^2(B',Z').
\end{equation}

Let $\sigma : G \to GL(V_\sigma)$ be an irreducible representation of
$G$ and $d_\sigma = \dim V_\sigma$, as before. Then
$$
        \sigma (Z') \subset C_{d_\sigma} = Z(GL(\CC,d_\sigma)) \cap
        SL(\CC,d_\sigma),
$$
which induces a morphism
$$
        \sigma_* : H^2(B', Z') \to H^2(B', C_{d_\sigma}) \simeq
        H^2(B', \ZZ/ d_\sigma \ZZ).
$$
{}From the above constructions, we obtain the following theorem.

\begin{teo}\ Let $\GR \to B$ be a bundle of compact, connected Lie groups on a
compact, connected smooth manifold $B$. Let $f_\sigma : B' \to B$ and
$Z' = G' \cap Z(G)$ be as before (Equations \eqref{eq.A} and
\eqref{eq.B}). The obstructions $\chi' \in H^2(B',Z')$ of Equation
\ref{eq.chi} and the Dixmier-Douady invariant $\chi_\sigma \in
H^2(B_\sigma, \ZZ/d_\sigma \ZZ)$ are related by $f_\sigma^*
(\chi_\sigma) = \sigma_* (\chi')$.
\end{teo}

\begin{proof}\
This follows from the definitions of the obstructions $\chi'$ and
$\chi_\sigma$, from the fact that the morphism $G \to GL(V_\sigma)$
maps $Z' \subset C_{d_\sigma}$, and from the naturality of the
boundary map in non-abelian cohomology.
\end{proof}\smallskip

The Dixmier-Douady invariant was recently been shown to be relevant
in the study of Ramond-Ramond fields (see \cite{FW, Mathai, W3, W4}
and the references therein). On the other hand, the algebra $C^*(\GR)$
is naturally a direct sum of algebras with controlled Dixmier-Douady
invariant. This suggests then the question whether Ramond-Ramond
fields can be obtained as indices of operators invariant with respect
to a bundle of compact Lie groups $\GR \to B$. The Dixmier-Douady
invariant of the resulting fields of algebras can be determined in
terms of a unique obstruction defined in terms of $\GR$, at least in
the case when the holonomy map $\pi_1(B,b_0) \to \Aut(G)/\Aut_0(G)$ is
trivial.

\section{The analytic index: an algebraic approach\label{S.An.Ind}}

We shall use below $\widehat{\otimes}$, the minimal tensor product of
$C^*$-algebras. Recall that the minimal tensor product of
$C^*$-algebras is defined to be (isomorphic to) the completion of the
image of $\pi_1 \otimes \pi_2$, the tensor product of two injective
representations $\pi_1$ and $\pi_2$. The same definition applies to
the tensor products $\widehat{\otimes}_C$ over a central subalgebra
$C$ \cite{Sakai}.

\begin{lem}\label{lemma.stable}\
Assume $\dim Y > \dim \GR$. Also, let ${\mathcal K} \to B$  denote
the locally trivial bundle of algebras whose fiber at $b \in B$ is
the algebra ${\mathcal K}(Y_b)$ of compact operators on $L^2(Y_b)$.
Define $C^*(\GR) \widehat\otimes \mathcal K =
C^*(\GR) \widehat{\otimes}_{C(B)} \Gamma( \mathcal K)$.
Then $C^*(Y,\GR)$ is Morita-equivalent to a direct summand of of $ C^*(\GR)
\widehat{\otimes} {\mathcal K}.$ Consequently, there is a natural
map
$$
        K_i(C^*(Y,\GR)) \to K_i( C^*(\GR) \widehat{\otimes} \mathcal K )
        \simeq K_i(C^*(\GR)).
$$
\end{lem}

\begin{proof}\
We shall regard $C^*(\GR) \widehat{\otimes} \mathcal K$ as an algebra
of operators acting on functions on $Y \times_B \GR$.  Let $\pi(g)$ be
the action of some $g \in \GR_b$ on $Y_b \times \GR_b$.  Denote by
$p_b = \int_{\GR_b} \pi(g) dg$, the integral being with respect to the
normalized Haar measure. Then $p_b ^2 = p_b$ is a self-adjoint
projection in the algebra $\mathcal B(L^2(Y_b)) \otimes C^*(\GR_b)$.
Let $p = (p_b)$. Then pointwise multiplication by $p_b$ defines a
multiplier of $C^*(\GR) \widehat{\otimes} \mathcal K$ (that is, it
maps $C^*(\GR) \widehat{\otimes} \mathcal K$ to itself by left or
right multiplication).

By a standard argument \cite{Rosenberg2}, the algebra $p(C^*(\GR)
\widehat{\otimes} \mathcal K)p$ is isomorphic to $C^*(Y,\GR)$. It is
then known that $p(C^*(\GR) \widehat{\otimes} \mathcal K)p$ is Morita
equivalent to $ (C^*(\GR) \widehat{\otimes} \mathcal K) p (C^*(\GR)
\widehat{\otimes} \mathcal K)$ \cite{Blackadar}. This completes the
proof.
\end{proof}\smallskip

We proceed now to define the index of an {\em elliptic, invariant family}
of operators
$$
        D = (D_b) \in M_N(\Psm{m}Y) = \Psm{m}{Y;\CC^N},
$$
without any holonomy assumption on $\GR \to B$.
We assume that $Y$ is compact, for simplicity; otherwise, we need
to use algebras with adjoint units. First, we observe that there exists an
exact sequence
\begin{equation}\label{eq.exact.seq}
  0 \to C^*(Y,\GR) \to {\mathcal E} \to \CI(S^*_{vt}Y) \to 0,\quad
  {\mathcal E}:= \Psm {0}Y + C^*(Y,\GR).
\end{equation}
The operator $D$ (or, rather, the family of operators $D = (D_b)$) has an
invertible principal symbol, and hence the family $T = (T_b)$,
$$
        T_b:=(1 + D_b^*D_b)^{-1/2}D_b,
$$
consists of elliptic, invariant operators, because the algebra of
pseudodifferential operators on a compact manifold is closed under holomorphic
functional calculus. Moreover, $T \in
{\mathcal E}=\Psm {0}Y + C^*(Y,\GR)$.  It's principal symbol is
still invertible, and hence defines a class $[T] \in
K_1(\CI(S^*_{vt}Y)) \simeq K^1(S^*_{vt}Y)$. Let
\begin{equation}\label{eq.boundary.map}
        \pa : K_1^{alg}(S^*_{vt}Y)) \to K_0^{alg}(C^*(Y,\GR))
        \simeq K_0 (C^*(Y,\GR))
\end{equation}
be the boundary map in the $K$-theory exact sequence
\begin{multline*}
        K_1^{alg}(C^*(Y,\GR)) \to K_1^{alg}({\mathcal E}) \to
        K_1^{alg}(S^*_{vt}Y)) \stackrel{\pa}{\rightarrow}K_0^{alg}(C^*(Y,\GR))
        \to K_0^{alg}({\mathcal E}) \to K_0^{alg}(S^*_{vt}Y))
\end{multline*}
associated to the exact sequence \ref{eq.exact.seq}. We hence obtain a
group morphism $\pa : K_1^{\alg}(\CI(S^*_{vt}Y)) \to
K_0(C^*(Y,\GR))$. By combining this morphism with the canonical
morphism $K_0(C^*(Y,\GR)) \to K_0(C^*(\GR))$, we obtain our desired
morphism,
\begin{equation}
        {\ind}_a : K_1^{\alg}(\CI(S^*_{vt}Y))  \to  K_0(C^*(\GR)),
\end{equation}
which we shall call {\em the analytic index morphism}. The image of
$D$ under the composition of the above maps, that is $\ind_a([T])$,
will be denoted $\ind_a(D)$ and will be called {\em the analytic index
of} $D$. The analytic index morphism descends to a group morphism
$K^1(S^*_{vt}Y) \to K_0(C^*(\GR))$ denoted in the same way. For $T$
acting between not necessarily equal vector bundles, we can proceed
similarly by using bivariant $K$-theory \cite{ConnesNCG, SkandExpo} to
define a morphism
\begin{equation}\label{eq.a.i.a}
        \ind_a : K^0(T^*_{vt}Y) \to K_0(C^*(\GR)).
\end{equation}

We still need to prove that the two definitions of the analytic index
coincide.

\begin{teo}\label{thm.equal}\ Let $\GR \to B$ be a bundle of
Lie groups with representation theoretic finite holonomy.  Then the
morphisms $a-\ind$ and $\ind_a$ defined in Equations \eqref{eq.a.i.g}
and \eqref{eq.a.i.a} are equal:
$$
        a-\ind=\ind_a : K^0(T^*_{vt}Y) \to K_0(C^*(\GR)).
$$
\end{teo}

\begin{proof}\
Let $F = (D_b)$ be a family in $\Psm{m}{Y;E_0,E_1}$. Choose $\KER$ as
in our first definition of the analytic index. Let $P$ be the
orthogonal projection onto the range of $\KER$. By a small
perturbation, we can arrange that $P \in \Psm{-\infty}{Y; \CC^N}$. We
can replace then $F$ by $F(1-P)$. A standard calculation (see
\cite{Blackadar} for example) then shows that $\pa [F] = [\KER] -
[\COK]$.
\end{proof}\smallskip

\section{Appendix}\label{banachcath}

Let us recall some general constructions of K-theory from \cite{KaKn}
and \cite{KarCli}. First we need some definitions.

\begin{dfn}\label{kar1}~\rm\cite[I.6.7]{KaKn}
An additive category $\cC$ is called
{\em pseudo-Abelian\/},
\label{ind:catpseud}
 if for each object $E$ from $\cC$ and each morphism $p:E\to E$,
satisfying to a condition $p^2=p$ (i.~e. an idempotent) there exists
the kernel $\Ker p$.  For an arbitrary additive category $\cC$ there
exists {\em associated pseudo-Abelian category\/}
\label{ind:catpseudass}
$\wt\cC$ which is a solution of the corresponding universal problem
and is defined as follows \cite[I.6.10]{KaKn}.  Objects of $\cC$ are
pairs $(E, p)$, where $E\in\Ob (\cC) $ and $p$ is a projector in
$E$. A morphism from $(E, p)$ to $(F, q)$ is such a morphism $f: E\to
F$ in $\cC$, that $f\circ p=q\circ f=f$.
\end{dfn}

\begin{dfn}
\label{kar2}~\rm\cite[\S~II.1]{KaKn}\label{ind:simmetrization}
Let $M$ be an abelian monoid.  On the product $M\times M$ we consider
the equivalence relation
$$
        (m,n)\sim
        (m',n')\:\Lra\:\exists\:p,q:\:(m,n)+(p,p)=(m',n')+(q,q).
$$
Let $S(M)$ be the quotient of $M \times M$ by the above equivalence
relation. Then $S(M)$ is a group, called the {\em group completion}.
If we consider now an additive category $\cC$ and denote by $\dot E$
the isomorphism class of an object $E$ from $\cC$, then the set $\F
(\cC) $ of these classes is equipped with a structure of an Abelian
monoid with respect to operation $\dot E + \dot F= (E\oplus
F)^\bullet$.  In this case the group $S (\F (\cC)) $ is denoted by $K
(\cC) $ and is called {\em Grothendieck group\/}
\label{ind:GrothGroup}
of the category $\cC$.
\end{dfn}

\begin{dfn}\label{kar3}~\rm\cite[\S~II.2] {KaKn}
A {\em Banach structure\/} on an additive category $\cC$ is defined by
a Banach space structure on all groups $\cC (E, F)$, where $E$ and $F$
are arbitrary objects from $\cC$ such that the composition of
morphisms $\cC (E, F) \times\cC (F, G) \to \cC (E, G) $ is bilinear
and continuous. We also say that $\cC$ is a {\em Banach category\/}.
\label{ind:catbanach}
\end{dfn}

\begin{dfn}\label{karpsab}~\rm\cite[\S~I.6]{KaKn}
Suppose that $\cC$ is an additive category. The category $\cC$ is
called {\em pseudo-abelian,\/} if, for each object $E$ of $\cC$ and
each morphism $p:E\to E$ satisfying the condition $p^2=p$, there
exists the kernel of $p$.
\end{dfn}

\begin{dfn}\label{kar4}~\rm\cite[\S~II.2]{KaKn}
Let $\cC$ and $\cC'$ be additive categories. An additive functor
$\f:\cC\to\cC'$ is called {\em quasi-surjective\/}
\label{ind:functquasi}
if each object of $\cC'$ is a direct summand of an object of type $\f
(E)$.  A functor $\f$ called {\em full\/}
\label{ind:functpoln}
if, for any $E,\,F\in\Ob (\cC)$, the map $\f (E, F) :\cC (E, F)
\to\cC' (\f (E),\f (F)) $ is surjective. For Banach categories $\f$ is
called {\em Banach\/}
\label{ind:functbanach}
if this map $\f (E, F)$ linear and continuous.
\end{dfn}

\begin{dfn}\label{kar5}~\rm\cite[II.2.13]{KaKn}
Let $\f:\cC\to\cC'$ be a quasi-surjective Banach functor. We shall
denote by $\G(\f)$ the set consisting of triples of the form $(E,
F,\a)$, where $E$ and $F$ are objects of the category $\cC$ and $\a:\f
(E) \to\f (F) $ is an isomorphism.  The triples $ (E, F,\a) $ and
$(E',F',\a') $ are called {\em isomorphic\/}, if there are
isomorphisms $f: E\to E'$ and $g: F\to F'$ such that the diagram
$$
        \xymatrix{ \f(E)\ar[r]^\a \ar[d]_{\f(f)}&
        \f(F)\ar[d]^{\f(g)}\\ \f(E')\ar[r]^{\a'}&\f(F')}
$$
commutes. A triple $(E,F,\a)$ is {\em elementary\/} if $E=F$ and
isomorphism $\a$ is homotopic in the set of automorphisms of $\f (E) $
to the identical isomorphism $\Id_{\f (E)}$.  We define {\em sum\/} of
two triples $ (E, F,\a) $ and $ (E ', F ',\a ')$ as
$$
(E\oplus E ', F\oplus F ',\a\oplus \a ').
$$
{\em The Grothendieck group $K (\f)$ of a functor $\f$\/} is defined
as quotient set of the monoid $\G(\f)$ with respect to the following
equivalence relation: $\s\sim\s'$ if and only if there exist
elementary triples $\t$ and $\t'$, that the triple $\s + \t$ is
isomorphic to the triple $\s' + \t'$.  The operation of addition
introduces on $K (\f) $ a structure of Abelian group. The class of a
triple we shall denote by $d (E, F,\a) $.
\end{dfn}

\begin{prop}\label{karopisnulfun}~\rm\cite[II.2.28]{KaKn}
Let $d(E,F,\a)$ be an element of $K(\f)$ where $\f:\cC\to\cC'$
is a full quasi-surjective Banach functor. Then $d(E,F,\a)=0$ if, and
only if, there exists an object $M$ of $\cC$ and an isomorphism
$\b: E\oplus M\to F\oplus M$ such that $\f(\b)=\a\oplus\Id_{\f(M)}$.
\end{prop}

\begin{dfn}\label{kar6}~\rm\cite[II.3.3]{KaKn}
Consider the set of pairs of the form $(E,\a)$, where $E$ is an object
of the category $\cC$ and $\a$ is an automorphism of $E$.  Two pairs
$(E,\a)$ and $(E',\a')$ are called {\em isomorphic\/}, if there is an
isomorphism $h: E\to E'$ in category $\cC$ such that the diagram
$$
        \xymatrix{ E\ar[r]^h \ar[d]_\a& E'\ar[d]^{\a'}\\ E\ar[r]^h&E'}
$$
commutes. The direct sum defines the operation of {\em addition\/} of
pairs. A pair $(E,\a)$ is called {\em elementary\/}, if the
automorphism $\a$ is homotopic to $\Id_E$ in the set of automorphisms
of $E$. Abelian group \cite[II.3.4]{KaKn} $K^{-1}(\cC)$ is defined as
a quotient set (with operation of addition) of the set of pairs
$\{(E,\a) \}$ with respect to the following equivalence relation:
$\s\sim\s'$ if and only if there are such elementary pairs $\t$ and
$\t'$, that $\s + \t$ is isomorphic to $\s' + \t'$.
\end{dfn}

\begin{dfn}\label{kar7}~\rm\cite[II.4.1] {KaKn}
Let $\cC$ be a Banach category and $C^{p, q}$ be the Clifford algebra.
We shall denote by $\cC^ {p, q}$ the category whose objects are pairs
$(E,\r)$, where $E\in\Ob (\cC) $ and $\r: C^ {p, q} \to\End (E) $ is a
homomorphism of algebras. A morphism from a pair $ (E,\r) $ to a pair
$ (E ',\r') $ is a $\cC$-morphism $f:E\to E'$ such that $f\circ\r
(\la)=\r(\la)\circ f$ for each element $\la\in\cC^{p, q} $.
\end{dfn}

Recall that there exist canonical morphisms $C^{p,q} \to C^{p,q+1}$.

\begin{dfn}\label{kar8}~\rm\cite[III.4.11]{KaKn}
Let $\cC$ be a pseudo-Abelian Banach category.  {\em The group\/}
$K^{p,q}(\cC)$ is defined as the Grothendieck group of the forgetful
functor $\cC^ {p, q + 1} \to\cC^ {p, q}$ (in the sense of the
Definition~\ref{kar5}).
\end{dfn}

The following statement can be easily obtained using the properties of
Clifford algebras.

\begin{teo}\label{kar9}~{\rm \cite[III.4.6, III.4.12]{KaKn}}
The groups $K^{p,q}(\cC)$ depend only on the difference
$p-q$. Besides, the groups $K^{0,0} (\cC)$ and $K^{0,1} (\cC) $ are
canonically isomorphic to groups $K (\cC)$ and $K^{-1}(\cC)$.
\end{teo}

\begin{dfn}\label{kar9a}
\rm Now we can define $K^{p-q} (\cC) =K^{p,q} (\cC)$ and similarly for
K-groups of functors.
\end{dfn}

We need also another description of K-groups, which is
equivalent~\cite[\S\S~III.4, III.5]{KaKn} to the initial.

\begin{dfn}\label{kar10}~\rm\cite [III.4.11, III.5.1]{KaKn}
Let $\cC$ be a pseudo-Abelian Banach category and let $E$ be a
$C^{p,q}$-module (an object of the category $\cC^{p,q}$).  A {\em
grading \/}
\label{ind:graduirovan}
of $E$ is an endomorphism
$\eta$ of  $E$ (considering as an object of $\cC$) such that
\begin{enumerate}
\item $\eta^2=1$,

\item $\eta\r (e_i) =-\r (e_i) \eta$, where $e_i$ are the generators
of Clifford algebra and $\r: C^{p,q} \to\End(E)$ is the homomorphism,
determining the $C^{p,q}$-structure on $E$.  \end{enumerate}
\end{dfn}

In other words, a grading of $E$ is a $C^{p,q+1}$-structure on $E$,
extending the initial $C^{p,q}$-structure (if we put $\r(e_{p + q +
1}) =\eta$).

Let us define the group $K^{p,q} (\cC) $ as the quotient group of the
free Abelian group, generated by the triples $(E,\eta_1,\eta_2)$,
where $E$ is a $C^{p,q}$-module and $\eta_1,\,\eta_2$ is a grading of
$E$ with respect to the subgroup, generated by relations
\begin{enumerate}
\item
$(E,\eta_1,\eta_2)\oplus(F,\xi_1,\xi_2)=(E\oplus
F,\eta_1\oplus\xi_1,\eta_2\oplus \xi_2) $,

\item $(E,\eta_1,\eta_2)=0$, if $\eta_1$ is homotopic to $\eta_2$ in
the set of gradations of $E$.
\end{enumerate}

As usual, by $d(E,\eta_1,\eta_2) \in K^{p, q} (\cC)$ we shall denote
the class of triple $ (E,\eta_1,\eta_2) $.


\providecommand{\bysame}{\leavevmode\hbox to3em{\hrulefill}\thinspace}

\end{document}